\documentclass[a4paper]{article}
\usepackage{RR}
\usepackage{hyperref}
\RRdate{Juillet 2009}

\usepackage{graphicx,a4wide}
\usepackage{times}
\usepackage{graphicx}
\usepackage{amsmath,amssymb,float}

\newcommand{\be} {\begin{equation}}
\newcommand{\ee} {\end{equation}}
\newcommand{\bea} {\begin{eqnarray}}
\newcommand{\eea} {\end{eqnarray}}
\newcommand{\Bea} {\begin{eqnarray*}}
\newcommand{\Eea} {\end{eqnarray*}}

\newcommand{\la} {\lambda}

\newcommand{\noi} {\noindent}

\newcommand{\var} {\varphi}
\newcommand{\R}{{\mathbb R}} 
\newcommand{\Z}{{\mathbb Z}} 
\newcommand{\cqfd}{\hfill \rule{1.5mm}{3mm}}
 \newtheorem{lemma}{Lemma}

\RRauthor{
Adimurthi\thanks[sfn]{TIFR-CAM, P.O. Box 1234, Bangalore 560012, India}
\thanks{E-mail: aditi@math.tifbng.res.in}
\and
G. D. Veerappa Gowda\thanksref{sfn}
\thanks{E-mail: gowda@math.tifbng.res.in}
\and
J\'er\^ome Jaffr\'e
\thanks{INRIA-Rocquencourt, B P 105, Le Chesnay Cedex, France,
E-mail: jerome.jaffre@inria.fr }}
\authorhead{Adimurthi \& G. D. Veerappa Gowda \& J. Jaffr\'e}
\RRtitle{Application du flux DFLU aux syst\`emes de lois de conservation} 
\RRetitle{Applications of the DFLU flux to systems of conservation laws}
\RRnote{This work was partially supported by the French-Indo cooperation project CEFIPRA 3401-2.}
\RRabstract{The DFLU numerical flux was introduced in order to solve hyperbolic scalar conservation laws with a flux function discontinuous in space. We show how this flux  can be used to solve systems of conservation laws. The obtained numerical flux is very close to a Godunov flux. As an example we consider a system modeling polymer flooding in oil reservoir engineering.
}
\RRresume{Le flux num\'erique DFLU a \'et\'e introduit afin de r\'esoudre des lois de conservation scalaires hyperbolique avec des fonctions de flux discontinues en espace. Nous montrons comment ce flux peut \^etre utilis\'e pour r\'esoudre des syst\`emes de lois de conservation. On obtient ainsi un flux num\'erique tr\`es proche du flux de Godunov. Comme exemple on consid\`ere un syst\`eme mod\'elisant l'injection de polym\`ere en ing\'eni\'erie de r\'eservoir p\'etrolier.}
\RRmotcle{Volumes finis, diff\'erences finies, solveurs de Riemann, syst\`emes de lois de conservation, \'ecoulements en milieu poreux, injection de polym\`eres.}
\RRkeyword{
Finite volumes, finite differences, Riemann solvers, system of conservation laws, flow in porous media, polymer flooding.
}
\RRprojets{Estime}
\RRtheme{\THNum}
\RCParis

\begin{document}
\RRNo{7009}
\makeRR

\section{Introduction}
The main difficulty in the numerical solution of systems of conservation laws is the complexity of constructing the Riemann solvers. One way to overcome this difficulty is to consider centered schemes as in \cite{LaxWen60,NesTad90,Toro99,Toro06,AdiGowJaf09}. However, in general these schemes are more diffusive than Godunov type methods based on exact or approximate Riemann solvers when this alternative is available. Therefore in this paper we will consider Godunov type methods. Most often the numerical solution requires the calculation of eigenvalues or eigenvectors of the Jacobian matrix of the system. This is even more complicated when the system  is non-strictly hyperbolic, i.e. eigenvectors are not linearly independent. In this paper we present a new approach which do not require such eigenvalue and eigenvector calculations.

Let us consider a system of conservation laws in conservative form
\[ \mathbf{U}_t + (\mathbf{F}(\mathbf{U}))_x=0, \quad \mathbf{U}=(u^1,\cdots,u^J), \quad \mathbf{F}=(f^1,\cdots,f^J). \]
A conservative finite volume method reads
\[ \dfrac{\mathbf{U^{n+1}_i}-\mathbf{U^n_i}}{\Delta t} + \dfrac{\mathbf{F}^{n}_{i+1/2} - \mathbf{F}^{n}_{i-1/2}}{h} =0 \]
where $\mathbf{F}^{n}_{i+1/2}$ is a numerical flux calculated using an exact or approximate Riemann solver. 
In a first order scheme this numerical flux is calculated using the left and right values 
$\mathbf{U}_i^n$ and $\mathbf{U}_{i+1}^n$. If we solve the equation field by field the $j$-th equation reads
\[ \dfrac{u^{j,n+1}_i-u^{j,n}_i}{\Delta t} + \dfrac{F^{j,n}_{i+1/2} - F^{j,n}_{i-1/2}}{h} =0\]
where the $j$-th numerical flux is a function of $\mathbf{U}^n_i$ and $\mathbf{U}^n_{i+1}$:
\[ F^{j,n}_{i+1/2} = F^j(u^{1,n}_i,\cdots, u^{j,n}_i,\cdots,u^{J,n}_i,u^{1,n}_{i+1},\cdots, u^{j,n}_{i+1},\cdots,u^{J,n}_{i+1}), \quad j=1,\cdots,J. \]
This flux function  can be calculated by solving the scalar Riemann problem for $t>t_n$:
\begin{align}
\label{dconsl}
&u^j_t + (\tilde{f}^{j,n}(u^j,x))_x = 0,\\ 
&u^j(x,t_n)=u^{j,n}_i \mbox{ if } x <x_{i+1/2}, \; u^j(x,t_n)=u^{j,n}_{i+1} \mbox{ if } x >x_{i+1/2}, \nonumber
\end{align}
where the flux function $\tilde{f}^j$, discontinuous at the point $x=x_{i+1/2}$, is defined by
\be \begin{array}{l} 
\tilde{f}^{j,n}(u^j,x) \equiv \tilde{f}^{j,n}_L(u^j) = f^j(u^{1,n}_i,\cdots, u^{j-1,n}_i,u^{j},u^{j+1,n}_i,\cdots,u^{J,n}_i) \mbox{ if } x <x_{i+1/2},\\
\tilde{f}^{j,n}(u^j,x) \equiv  \tilde{f}^{j,n}_R(u^j) = f^j(u^{1,n}_{i+1},\cdots, u^{j-1,n}_{i+1},u^{j},u^{j+1,n}_{i+1},\cdots,u^{J,n}_{i+1}) \mbox{ if } x >x_{i+1/2}
\label{fluj}
  \end{array} \ee
(L and R refer to left and right of the point $x_{i+1/2}$).

Scalar conservation laws like equation (\ref{dconsl}) with a flux function discontinuous in space have been the object of many studies \cite{chacohjaf87,mochen87,lantvewin92,gimris92,jaf96,kaa99,Tow00,Tow01,BurKarRisTow03,SegVov03,AdiJafGow04, SidJaf09}. In particular, in \cite{AdiJafGow04} a Godunov type finite volume scheme was proposed and convergence to a proper entropy condition was proved, provided that the left and right flux functions have exactly one local maximum and the same end points (the case where the flux functions has exactly one local minimum can be treated by symetry). At the discontinuity the interface flux, that we call the DFLU flux, is given by the formula 
\be F^n_{i+1/2}(u_L,u_R) =  \min\Bigl\{ f_L(\min\{u_L,\theta_L\}), f_R(\max\{u_R,\theta_R\})\Bigr\}, 
\label{dflu} \ee
 if $f$ denotes the scalar flux function and $\theta_L=$argmax$(f_L)$, $\theta_R=$argmax$(f_R)$. When $f_L \equiv f_R$ this formula is equivalent to the Godunov flux so formula (\ref{dflu}) can be seen as an extension of the Godunov flux to the case of a flux function discontinuous in space.
In the case of systems formula (\ref{dflu}) can be applied to the fluxes $\tilde{f}^{j,n}_L$ and 
$\tilde{f}^{j,n}_R$.

To illustrate the method we consider the system of conservation laws arising for polymer flooding in reservoir simulation which is described in section \ref{polymer}.
This system, or similar systems of equations, is nonstrictly hyperbolic and is studied in several papers \cite{Tem82, JohWin88, JohTveWin89,IssTem90}. For example in \cite{JohWin88} the authors solve Riemann problems associated to this system when gravity is neglected and therefore the fractional flow function is an increasing function of the unknown.  In this case, the eigenvalues of the corresponding Jacobian matrix  are positive and hence it is less difficult to construct Godunov type schemes which turn out to be upwind schemes. When  the above model with gravity effects is considered, then the flux function is not necessarily monotone and hence the eigenvalues can change sign. This makes the construction of Godunov type schemes more difficult as it involves exact solutions of Riemann problems with a nonmonotonous fractional flow function.
Therefore in section \ref{Riemann} we solve the  Riemann problems in the general case when gravity terms are taken into account so the flux function is not anymore monotone. This will allow to compare our method with that using an exact Riemann solver.
In section \ref{finite difference} we consider Godunov type finite volume schemes. We present the DFLU scheme for the system of polymer flooding and compare it to the Godunov scheme whose flux is given by the exact solution of the Riemann problem. We also present several other possible numerical fluxes, centered like Lax-Friedrichs or Force, or upstream like the upstream mobility flux commonly used in reservoir engineering \cite{AziSet79, BreJaf91,SidJaf09}.
Finally in section \ref{numericalresults} we compare numerically the DFLU method with these fluxes.

\section{A system of conservation laws modeling polymer flooding}
\label{polymer}
A polymer flooding model for enhanced oil recovery in petroleum engineering was introduced in \cite{Pop80} as the following $2\times 2$ system of conservation laws
\be
\begin{array}{rrll}
 s_t + f(s, c)_x & =& 0 \\
 (sc+a(c))_t+(cf(s,c))_x &=&0
 \label{polyeq1}
\end{array}
\ee
where $t > 0$ and $x \in \R$, $(s,c) \in I \times I$ with $I=[0,1]$. $s=s(x,t)$ denotes the saturation of the wetting phase, so $1-s$ is  the saturation of the oil phase. $c=c(x,t)$ denotes the concentration of the polymer in the wetting phase which we have normalized. Here the porosity was set to 1 to simplify notations. The flux function $f$ is the Darcy velocity of the wetting phase  $\var_1$ and is determined by the relative permeabilities and the mobilities of the wetting and oil phases, and by the influence of gravity:
\be
\label{fg-twophase}
f(s,c)=  \var_1 = \dfrac{\lambda_1(s,c)}{\lambda_1(s,c) + \lambda_2(s,c)} [ \var + (g_1-g_2)\lambda_2(s,c) ].
\ee
The quantities $\lambda_\ell, \ell=1,2$  are the mobilities of the two phases, with $\ell =1$ referring to the wetting phase and  $\ell =2$ referring to the oil phase:
\[ \lambda_\ell(s,c) = \dfrac{K kr_\ell(s)}{\mu_l(c)}, \ell=1,2, \]
where $K$ is the absolute permeability, and $kr_\ell$ and $\mu_\ell$ are respectively the relative permeability and the viscosity of the phase $\ell$.
$kr_1$ is an increasing function of $s$ such that $kr_1(0)=0$ while $kr_2$ is a decreasing function of $s$ such that $kr_2(1)=0$. Therefore $\lambda_\ell, \ell=1,2$ satisfy
\begin{equation}
\label{eq3}
\begin{array}{l}
\lambda_1=\lambda_1(s,c) \mbox{is an  increasing functions of} \, s ,\;
\lambda_1(0,c) = 0\;  \forall c \in [0,1],\\
\lambda_2= \lambda_2(s,c) \, \mbox{is a decreasing functions of} \, s ,\;
\lambda_2(1,c)= 0\;  \forall c \in [0,1].
\end{array}
\end{equation}
The idea of polymer flooding is to dissolve a polymer in the injected water in order to increase the viscosity of the injected wetting phase. Thus the injected wetting phase will not be able to bypass oil so one obtains a better displacement of the oil by the injected phase. Therefore $\mu_1(c)$ is increasing with $c$ while $\mu_2$ will be taken as a constant assuming there is no chemical reaction between the polymer and the oil.  Therefore $f$ will decrease with respect to $c$.
The function $a=a(c)$ models the adsorption of the polymer by the rock and is increasing with $c$. 

$\var$ is the total Darcy velocity, that is the sum of the Darcy velocities of the
two phases $\var_1$ and $\var_2$:
\[ \var = \var_1+\var_2, \quad \var_1 = \dfrac{\lambda_1}{\lambda_1 + \lambda_2} [ \var + (g_1-g_2)\lambda_2], \quad
\var_2 = \dfrac{\lambda_2}{\lambda_1 + \lambda_2} [ \var + (g_2-g_1)\lambda_1].\]
$\var $ is a constant in space since we assume that the flow is incompressible.
The gravity constants $g_1, g_2$ of the phases are proportional to their density.

To equation (\ref{polyeq1})  we add the initial condition
\begin{equation}
(s(x,0), c(x, 0)) = (s_0 (x),c_0(x)). \label{inicond}
\end{equation}

\begin{figure}[htbp]
\begin{center}
\begin{minipage}[c]{5cm}
\setlength{\unitlength}{0.7pt}
\begin{picture}(100,110)(0,0)
\put(0,0){\includegraphics[height=3cm]{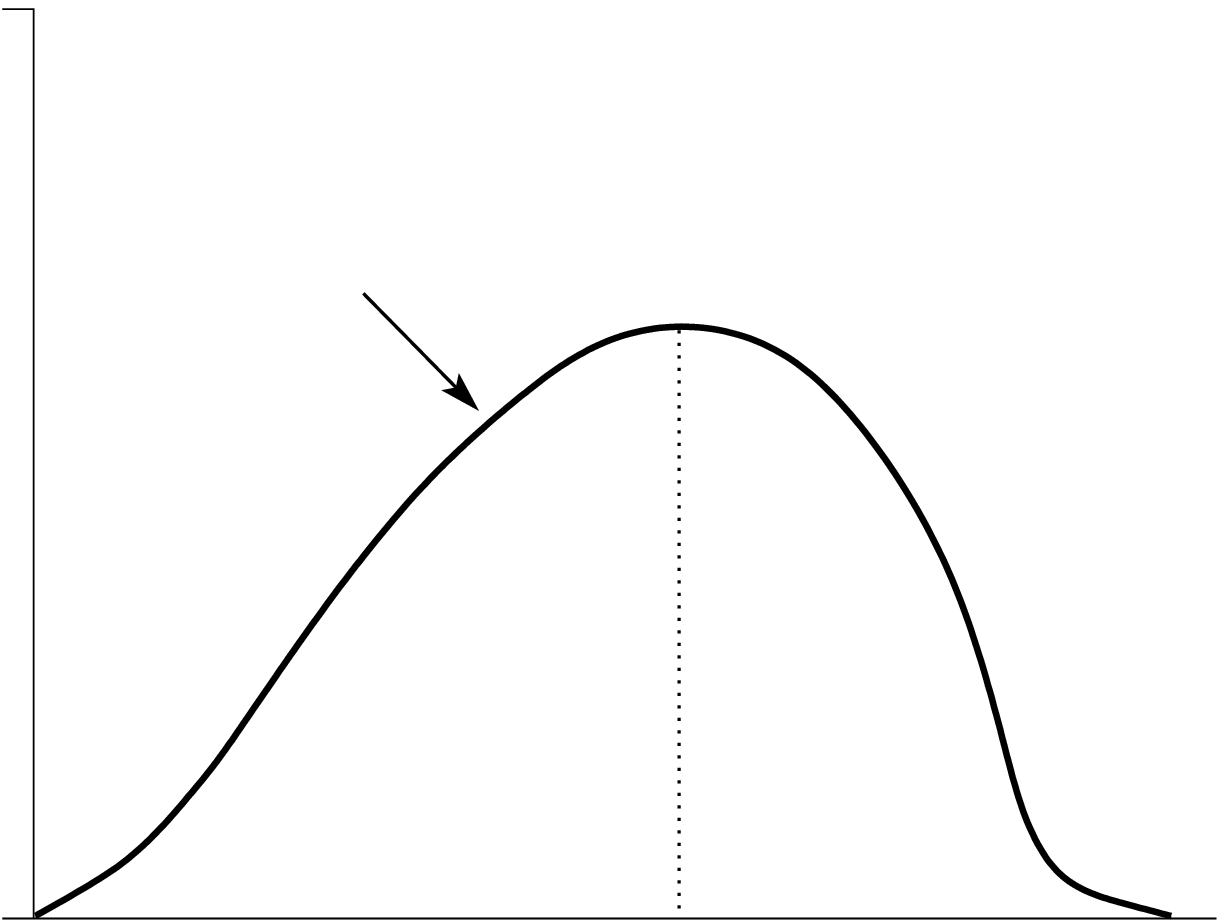}}
\put(45,90){\makebox(0,0){$f(\cdot,c)$}}
\put(-5,5){\makebox(0,0){0}}
\put(8,-7){\makebox(0,0){0}}
\put(91,-7){\makebox(0,0){$\theta$}}
\put(155,-7){\makebox(0,0){$1$}}
\put(164,-8){\makebox(0,0){$s$}}
\end{picture}
\end{minipage}
\begin{minipage}[c]{4.5cm}
\setlength{\unitlength}{0.7pt}
\begin{picture}(100,110)(0,0)
\put(0,0){\includegraphics[height=3cm]
{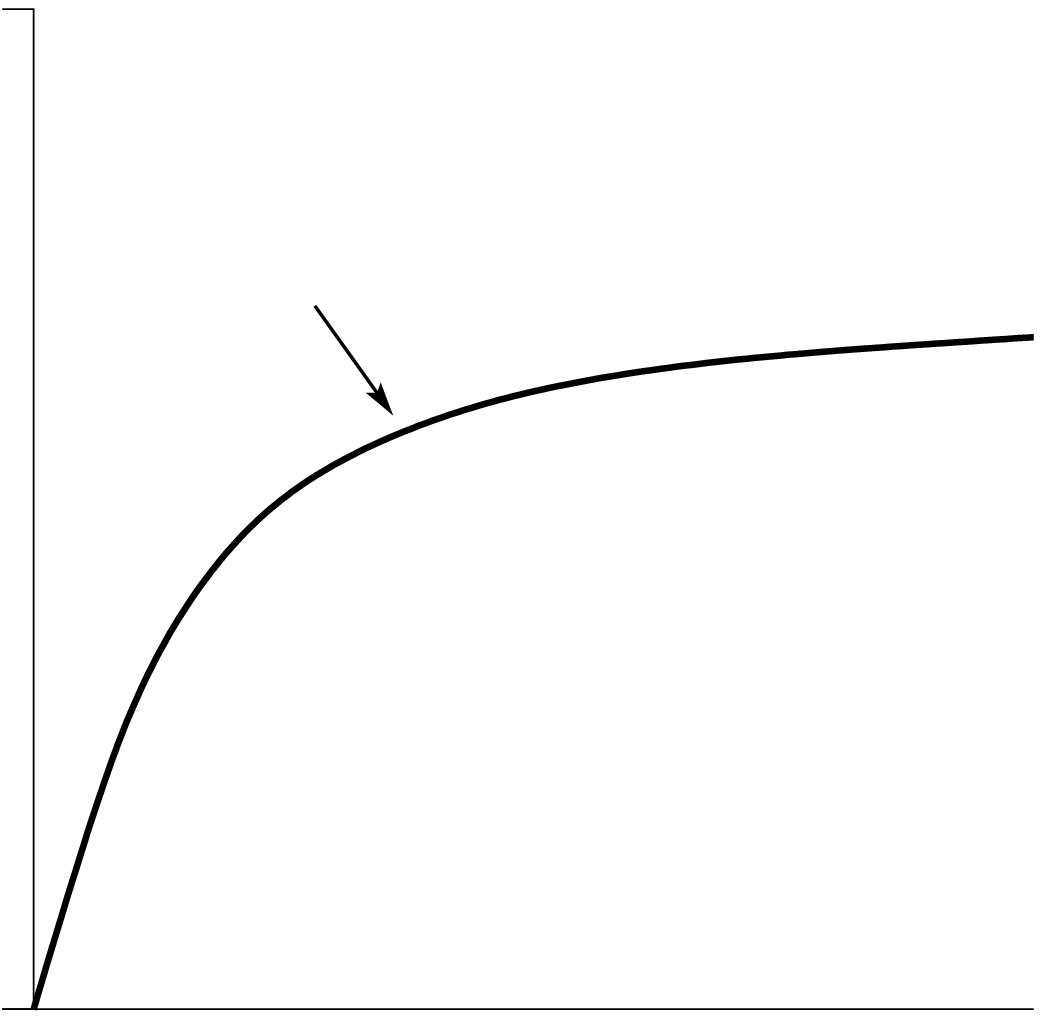}}
\put(35,95){\makebox(0,0){$a$}}
\put(-5,7){\makebox(0,0){0}}
\put(5,-8){\makebox(0,0){0}}
\put(127,-8){\makebox(0,0){$c$}}
\end{picture}
\end{minipage}
\caption{Shapes of flux function $s \rightarrow f(s,c)$ (left) and adsorption function $c \rightarrow a(c)$ (right).}
\label{shape-f-a}
\end{center}
\end{figure}

Since the case when $f$ is monotone was already studied in \cite{JohWin88,JohTveWin89}, we concentrate on the nonmonotone case which is more complicated and corresponds to taking into account gravity. Therefore we assume that $\var = 0$ so for the
nonlinearities of the system (\ref{polyeq1}). We will assume also that phase 1 is heavier than phase 2 ($g_1>g_2$) so  we can assume the following properties:
\begin{enumerate}
\item[(i)] $f(s,c) \ge 0, f(0,c)= f(1,c)=0$ for all $c \in I$.
\item[(ii)] The function $s \rightarrow f(s,c)$ has exactly one global
 maximum in $I$ with $\theta=$argmax$(f)$.
\item[(iii)] $f_c(s,c) < 0  \,\,\forall \,\,s \in (0,1)$ and for all $ c \in I$
\item[(iv)] The adsorbtion term $a=a(c)$ satisfies\\
   $a(0)=0,\,\,\,\, h(c)=\dfrac{da}{dc}(c)>0, \quad \dfrac{d^2a}{dc^2}(c) < 0$
 for all $c \in I$.  
 \end{enumerate}
 Typical shapes of functions $f$ and $a$ are shown in Fig. \ref{shape-f-a}.
We expand the derivatives in equations  (\ref{polyeq1}) and we plug the resulting first equation into the second one. Then we obtain the system in nonconservative form
 \[
\begin{array}{rrll}
 s_t + f_s(s, c) s_x + f_c(s,c) c_x & =& 0, \\
 (s +a'(c)) c_t+f(s,c)) c_x &=&0.
\end{array}
\]
 Let $U$ denote the state vector $U=(s,c)$ and introduce the upper triangular matrix
  \[  A(U) =  \left( \begin{array}{cc}
  f_s & f_c \\[0.2cm] 0 & \dfrac{f}{s + a'(c)}\end{array} \right) \]
and the system (\ref{polyeq1}) can be read in matrix form as
$$ U_t+ A(U)\,U_x\,=\,0.$$   
  
 The eigenvalues of $A$ are $\lambda^s=f_s$ and $\lambda^c=\dfrac{f}{s+a'}$, with corresponding eigenvectors $e^s=(1,0), e^c=(f_c,\lambda^c-\lambda^s) $ if $0<s<1$   and $e^c=(0,1)$ if $s=0,1$. The
 eigenvalue  $\lambda^s$ may change sign whereas the eigenvalue $\lambda^c$ 
 is always positive. One can observe that for each $c \in I$ there exists a unique $s^*=s^*(c) \in (0,1)$ such that
$$ \lambda^c(s^*,c)=\lambda^s(s^*,c)$$
(see Fig.\ref{rs0}). For this couple $(s^*,c)$, $\lambda^c=\lambda^s$, hence eigenvectors are not linearly independent and the problem is nonstrictly hyperbolic.

 Any weak solution of (\ref{polyeq1}) has to satisfy the Rankine-Hugoniot
 jump  conditions given by
\be
\begin{array}{rrll}
 f(s_R,c_R)-f(s_L,c_L) & =& \sigma(s_R -s_L),\\
 c_R f(s_R,c_R)-c_L f(s_L,c_L)&=&\sigma (s_R c_R +a(c_R) -s_L c_L - a(c_L)),
\label{RHcondition1}
\end{array}
\ee
where $(s_L,c_L), (s_R,c_R)$ denote the left and right values of the couple $(s,c)$ at a certain point of discontinuity.

When $c_R=c_L$, the second equation reduces to the first equation and the speed of the discontinuity $\sigma$ is given by the first equation only. Now we are interested in the case $c_R \neq c_L$.
By combining the two equations (\ref{RHcondition1}) we may write
$$ (c_R - c_L)f(s_L,c_L)=\sigma(c_R - c_L)s_L + \sigma( a(c_R) - a(c_L)) $$
where
\[ \sigma=\frac{f(s_L,c_L)}{s_L+\bar{a}_L(c_R)}, \quad \bar{a}_L(c)=\left \{\begin{array}{lll}
\dfrac{a(c)-a(c_L)}{c-c_L} &\mbox{if} &c \neq c_L,\\
a'(c) &\mbox{if}& c=c_L. \end{array} \right. 
\]
Plugging this into first equation of (\ref{RHcondition1}), we obtain
$$\sigma(s_R + \bar{a}_L(c_R))=\sigma(s_L + \bar{a}_L(c_R))+f(s_R,c_R)-f(s_L,c_L)=f(s_R,c_R).$$
Hence when $c_L \neq c_R$ the Rankine-Hugoniot condition ({\ref{RHcondition1}})
reduces to
\begin{equation}
\dfrac{f(s_R,c_R)}{s_R+\bar{a}_L(c_R)}=\dfrac{f(s_L,c_L)}{s_L +\bar{a}_L(c_R)}=\sigma.
\label{rhcondition2}
\end{equation}

\section{Riemann problem}
\label{Riemann}
In this section we  solve the Riemann problems associated with our system, that we solve system (\ref{polyeq1}) with the initial condition
\be s(x,0) = \left\{ \begin{array}{lll}
s_L &\mbox{if}& x<0, \\ s_R &\mbox{if}& x>0 \end{array} \right.  , \quad
c(x,0) = \left\{ \begin{array}{lll}
c_L &\mbox{if}& x<0, \\ c_R &\mbox{if}& x>0 \end{array} \right.  .
\label{riemannpb}
\ee
Solution to (\ref{riemannpb}) is constructed by using elementary waves associated with the system. There are two families of waves, refered to as the $s$ and $c$ families. $s$ waves consist of rarefaction and shocks (or contact discontinuity) across which $s$ changes continuously and discontinuously respectively, but across which $c$ remains constant. $c$ waves consist solely of contact discontinuities, across which both $s$ and $c$ changes such that $\dfrac{f(s,c)}{s+a'(c)}$ remains constant in the sense of (\ref{rhcondition2}).

First define a function $\bar{a}_L$ by
\[ \bar{a}_L(c) = \left \{\begin{array}{lll}
\dfrac{a(c)-a(c_L)}{c-c_L} &\mbox{if} &c \neq c_L,\\
a'(c) &\mbox{if}& c=c_L. \end{array} \right.
\]
                                                                               
We will restrict to the case $c_L > c_R$. The case  $c_L <  c_R$ can be treated similarly.

\begin{figure}[H]
\begin{center}
 \begin{picture}(200,100)(0,-10)
 \put(0,0){\includegraphics[height=3cm]{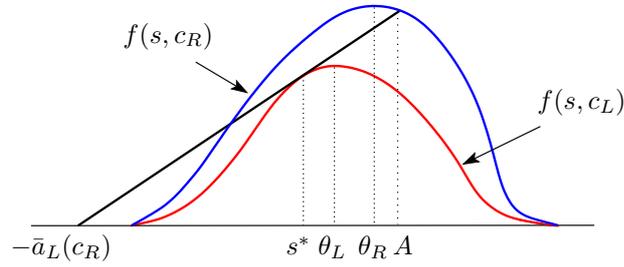}}
 \put(35,70){$f(s,c_R)$}
 \put(190,45){$f(s,c_L)$}
 \put(-7,-10){$-\bar{a}_L(c_R)$}
 \put(96,-10){$s^*$}
  \put(123,-10){$\theta_R$}
   \put(108,-10){$\theta_L$}
   \put(136,-10){$A$}
 \end{picture}
 \end{center}
 \caption{Two flux functions $f(s,c_L)$ and  $f(s,c_R)$ with $c_L > c_R$.}
 \label{rs0}
 \end{figure}

When  $c_L > c_R$ the flux functions for the first equation (\ref{polyeq1}) $s \rightarrow f(s,c_L)$ and $s \rightarrow f(s,c_R)$ are as represented in 
Fig. \ref{rs0}, that is $f(s,c_L) \leq f(s,c_R) \; \forall s \in (0,1)$. Let 
$\theta_L$ and $\theta_R$ be the points at which $f(.,c_L)$ and $f(.,c_R)$ reach their
maxima respectively.

 Let $s^* \in (0,1)$ be a point at which  
$ f_s(s^*,c_L)=\dfrac{f(s^*,c_L)}{s^*+\bar{a}_L(c_R)}$. Now draw a line through the points
$(-\bar{a}_L(c_R),0)$ and $(s^*,f(s^*,c_L)$ which intersects the curve $f(s,c_R)$ at a point $A \geq s^*$ (see Fig. \ref{rs0}).

 Our study of Riemann problems separates into two cases $s_L < s^* $ and $s_L \geq s^*$ which themselves separate into several subcases.

 \begin{itemize}
 \item {\bf Case 1:} $ s_L < s^*$.

    Draw a line passing through the points $(s_L,f(s_L,c_L))$ and $(-\bar{a}_L(c_R),0)$. This line intersects the curve $f(s,c_R)$ at points $\overline{s}$ and $B$ (see Fig. \ref{figR1a} ). Now we divide this into two subcases.  
 
 \item {\bf Case 1a:} $s_R < B$\\
  (a)  Connect $(s_L,c_L)$ to $(\overline{s},c_R)$ by $c$-wave with a speed 
 $$\sigma_c=\frac{f(s_L,c_L)}{s_L+\bar{a}_L(c_R)}=\frac{f(\overline{s},c_R)}{\overline{s}+\bar{a}_L(c_R)}. $$
  (b) Next connect  $(\overline{s},c_R)$ to $(s_R,c_R)$ by a $s$-wave, along the curve
     $f(s,c_R)$ (see Fig. \ref{figR1a}).

   For example if $ s_R \geq \overline{s}$ and $f(s,c_L)$ and $f(s,c_R)$ are concave functions then the solution of the Riemann problem is given by
\be (s(x,t),c(x,t)) = \left\{ \begin{array}{lll}
(s_L,c_L) &\mbox{if}& x<\sigma_c t, \\ (\overline{s},c_R) &\mbox{if}& \sigma_c t < x < \sigma_s t, \\ (s_R,c_R) &\mbox{if}& x > \sigma_s t, \end{array} \right. 
\ee
where
 $$\sigma_c=\frac{f(s_L,c_L)}{s_L+\bar{a}_L(c_R)}=\frac{f(\overline{s},c_R)}{\overline{s}+\bar{a}_L(c_R)},\quad \sigma_s=\frac{f(\overline{s},c_R) - f(s_R,c_R)}{\overline{s}-s_R}.$$

\noi Note that $ 0 < \sigma_c < \sigma _s$.     
\begin{figure}[H]
 \begin{picture}(320,100)(-40,-10)
 \put(0,0){\includegraphics[height=3cm]{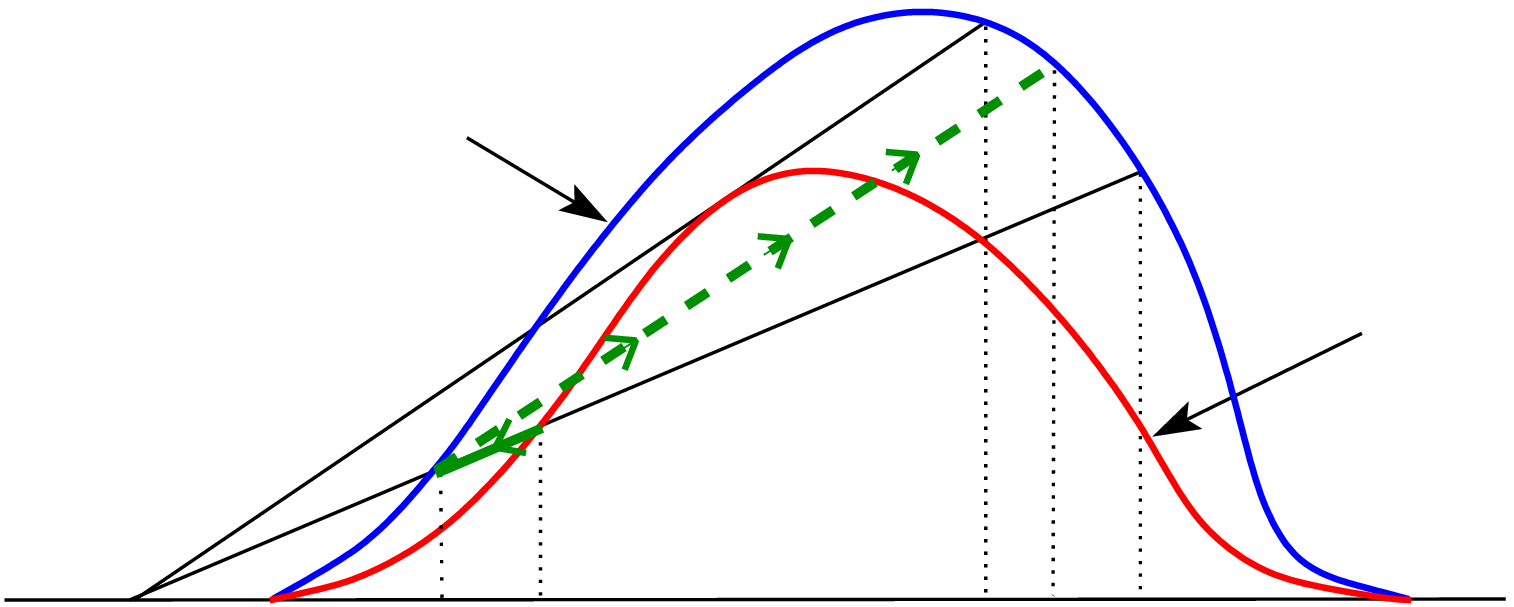}}
 \put(40,72){$f(s,c_R)$}
\put(96,-10){$s^*$}
 \put(190,45){$f(s,c_L)$}
 \put(60,-10){$\overline{s}$}
 \put(70,-10){$s_L$}
 \put(132,-10){$A$}
 \put(144,-10){$s_R$}
 \put(157,-10){$B$}

 \put(-7,-10){$-\bar{a}_L(c_R)$}
 \put(220,0){\includegraphics[height=1.8cm]{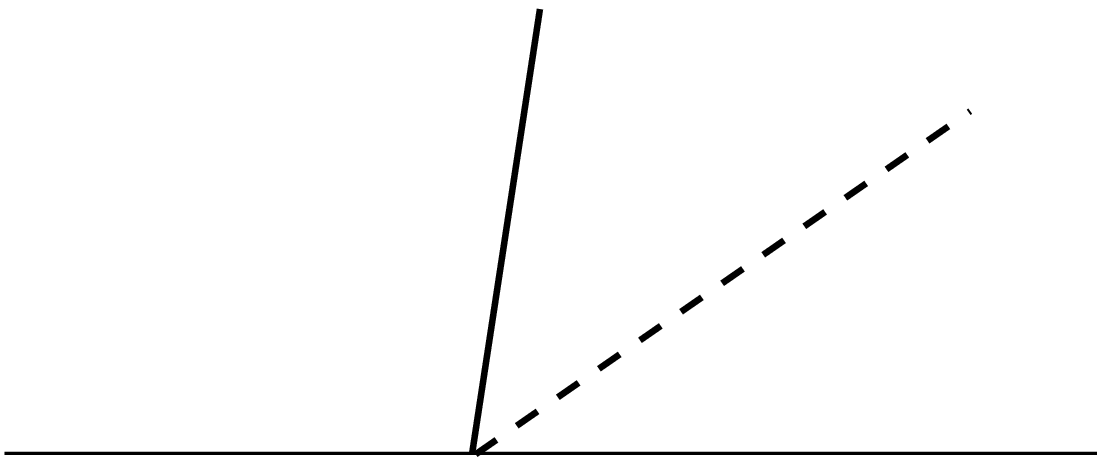}}
  \put(232,15){$(s_L,c_L)$}
 \put(228,-10){$(s_L,c_L)$}
 \put(271,-10){0}
  \put(300,-10){$(s_R,c_R)$}
  \put(310,13){$(s_R,c_R)$}
  \put(283,35){$(\overline{s},c_R)$}
  \put(255,57){$x=\sigma_c t$}
   \put(320,44){$x=\sigma_s t$}
 \end{picture}
 \caption{Solution of Riemann problem (\ref{riemannpb}) with $s_L <  s^*$ and $ s_R <  B$.}
 \label{figR1a}
 \end{figure}
 
 \item {\bf Case 1b:} $s_R \geq B $.

  Draw a line passing through the points $(s_R,f(s_R,c_R))$ and $(-\bar{a}_L(c_R),0)$. This line intersects the curve $f(s,c_L)$ at a point $\overline{s}$ (see Fig. \ref{rs1b}).

 (a)  Connect $(s_L,c_L)$ to $(\overline{s},c_L)$ by a $s$-wave along the curve $f(s,c_L)$.\\
 (b) Next connect  $(\overline{s},c_L)$ to $(s_R,c_R)$ by a $c$-wave with a speed $$\sigma_c=\frac{f(s_R,c_R)}{s_R+\bar{a}_L(c_R)}=\frac{f(\overline{s},c_L)}{\overline{s}+\bar{a}_L(c_R)}. $$
   
For example if $f(s,c_L)$ and $f(s,c_R)$ are concave functions then the solution is given by
\be (s(x,t),c(x,t)) = \left\{ \begin{array}{lll}
(s_L,c_L) &\mbox{if}& x<\sigma_s t, \\ (\overline{s},c_L) &\mbox{if}& \sigma_s
t < x < \sigma_c t \\ (s_R,c_R) &\mbox{if}& x > \sigma_c t \end{array} \right.
\ee
where
 $$\sigma_c=\dfrac{f(s_R,c_R)}{s_R+\bar{a}_L(c_R)}=\dfrac{f(\overline{s},c_L)}{\overline{s}+\bar{a}_L(c_R)}, \quad \sigma_s=\frac{f(\overline{s},c_L) - f(s_L,c_L)}{\overline{s}-s_L}.$$
\noi Note that $  \sigma_s < \sigma _c$ and $(s_L,c_L)$ is connected to
$(\overline{s},c_L)$ by a $s$-shock wave and $(\overline{s},c_L)$ is connected
to $(s_R,c_R)$ by a $c$-shock wave.

\begin{figure}[H]
 \begin{picture}(320,100)(-40,-10)
 \put(0,0){\includegraphics[height=3cm]{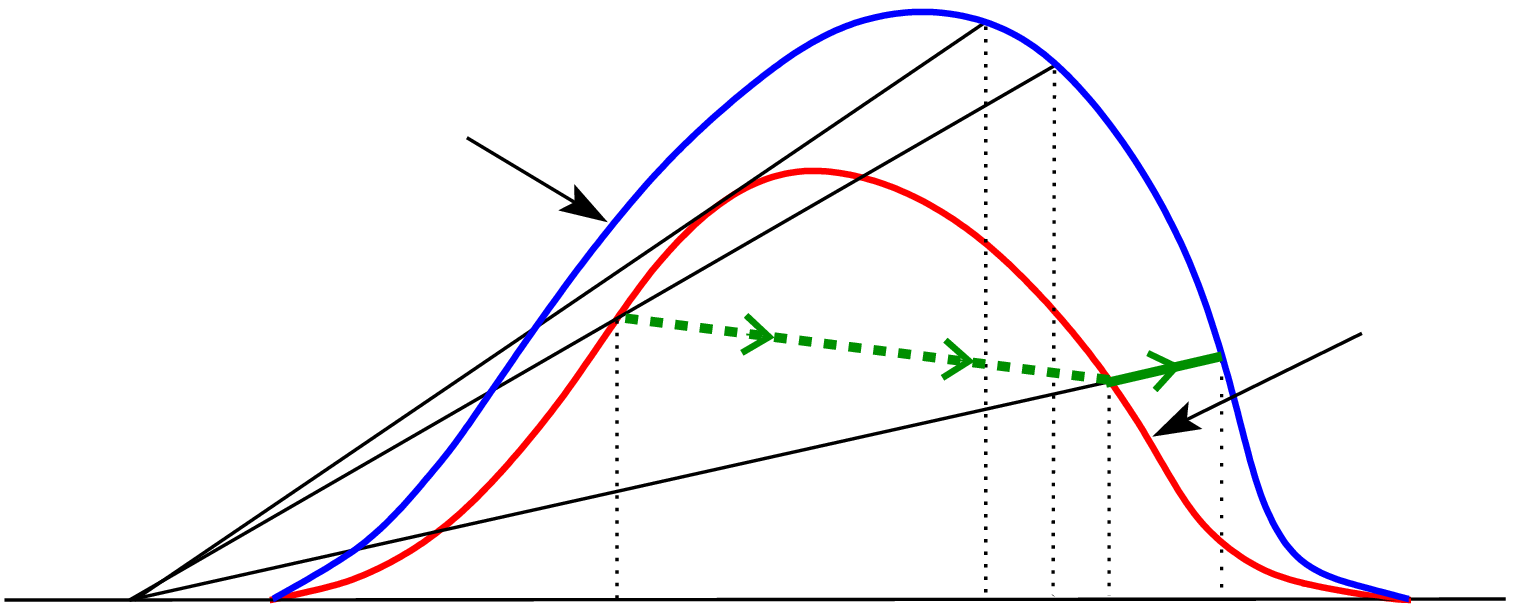}}
 \put(40,72){$f(s,c_R)$}
 \put(190,45){$f(s,c_L)$}
 \put(83,-10){$s_L$}
\put(96,-10){$s^*$}
 \put(132,-10){$A$}
 \put(144,-10){$B$}
 \put(157,-10){$\overline{s}$}
 \put(165,-10){$s_R$}
 \put(-7,-10){$-\bar{a}_L(c_R)$}
 \put(220,0){\includegraphics[height=1.8cm]{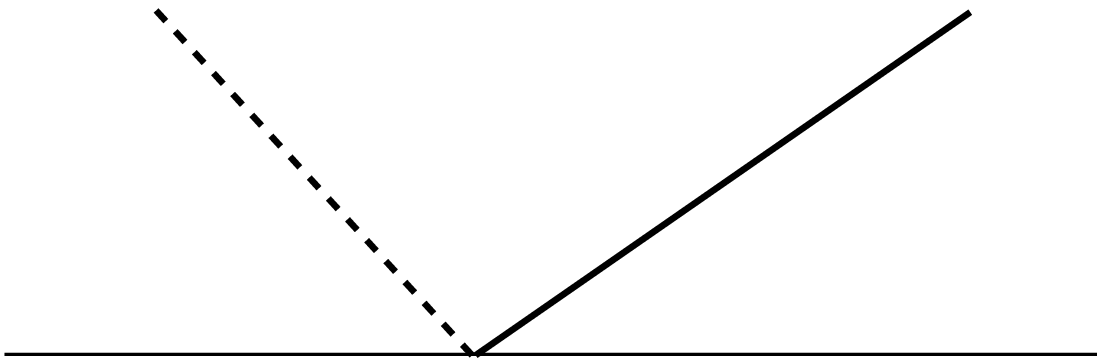}}
 \put(225,15){$(s_L,c_L)$}
 \put(230,-10){$(s_L,c_L)$}
 \put(285,-10){0}
  \put(315,-10){$(s_R,c_R)$}
  \put(325,15){$(s_R,c_R)$}
  \put(283,28){$(\overline{s},c_L)$}
  \put(243,55){$x=\sigma_s t$}
   \put(340,55){$x=\sigma_c t$}
 \end{picture}
 \caption{Solution of Riemann problem (\ref{riemannpb}) with $s_L <  s^*$ and
 $s_R \geq B $.}
 \label{rs1b}
 \end{figure} 
\end{itemize} 
\begin{itemize}
 \item {\bf Case 2:} $ s_L \geq s*$.\\

  \item {\bf Case 2a:} $s_R \leq A$ .\\
(a) Connect $(s_L,c_L)$ to $(s^*,c_L)$ by a $s$-wave along the curve $f(s,c_L)$.\\
(b) Connect $(s^*,c_L)$ to $(\overline{s},c_R)$ by a $c$-wave.\\ 
(c) Connect $(\overline{s},c_R)$ to $(s_R,c_R)$ by a $s$-wave along  the curve $f(s,c_R)$ (see Fig. \ref{rs2a}).

    For example if $s_R \leq \overline{s}$, then the solution is given by 

\[ (s(x,t),c(x,t)) = \left\{ \begin{array}{lll}
(s_L,c_L) &\mbox{if}& x<\sigma_1 t, \\ ((f_s)^{-1}(\frac{x}{t},c_L),c_L) &\mbox{if}& \sigma_1 t < x < \sigma_2 t ,\\ (\overline{s},c_R) &\mbox{if}& \sigma_2 t <x < \sigma_3 t,\\
((f_s)^{-1}(\frac{x}{t},c_R),c_R) &\mbox{if}&  \sigma_3 t < x < \sigma_4 t,\\
(s_R,c_R) &\mbox{if}&  x > \sigma_4 t,
 \end{array} \right.
\]
where
\[ \sigma_1=f_s(s_L,c_L), \quad \sigma_2 = f_s(s^*,c_L)=\frac{f(s^*,c_L)}{s^*+\bar{a}_L(c_R)}, \quad
 \sigma_3 = f_s(\overline{s},c_R), \quad \sigma_4= f_s(s_R,c_R).\]

   Here $ (s_L,c_L)$ is connected to $(s^*,c_L)$ by a $s$-rarefaction wave, $(s^*,c_L)$ is connected to $(\overline{s},c_R)$ by a $c$-shock wave and  (see Fig. \ref{rs2a}). If $s_R > \overline{s}$ then 
 $(\overline{s},c_R)$ would be connected to $(s_R,c_R)$ by a $s$-chock wave.
\begin{figure}[H]
 \begin{picture}(320,100)(-40,-10)
 \put(0,0){\includegraphics[height=3cm]{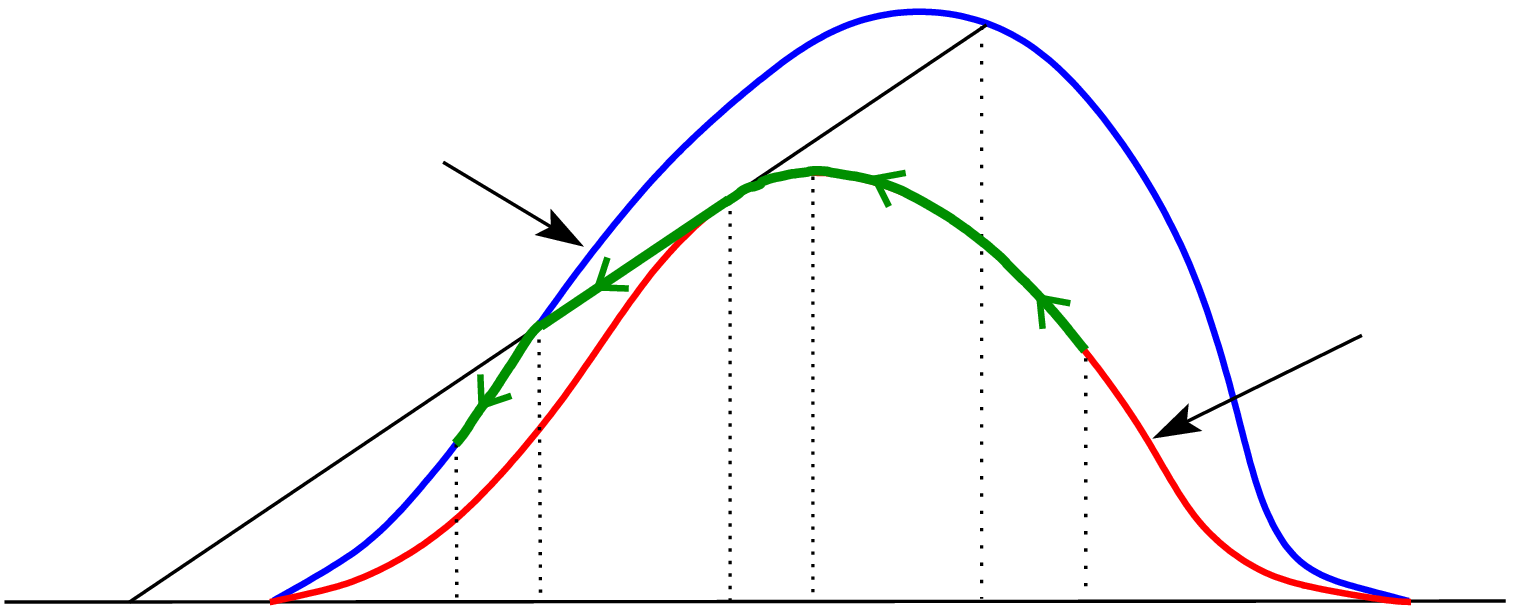}}
 \put(35,70){$f(s,c_R)$}
 \put(190,45){$f(s,c_L)$}
 \put(-8,-10){$-\bar{a}_L(c_R)$}
 \put(58,-10){$s_R$}
  \put(75,-10){$\overline{s}$}
 \put(96,-10){$s^*$}
   \put(111,-10){$\theta_L$}
   \put(220,0){\includegraphics[height=1.8cm]{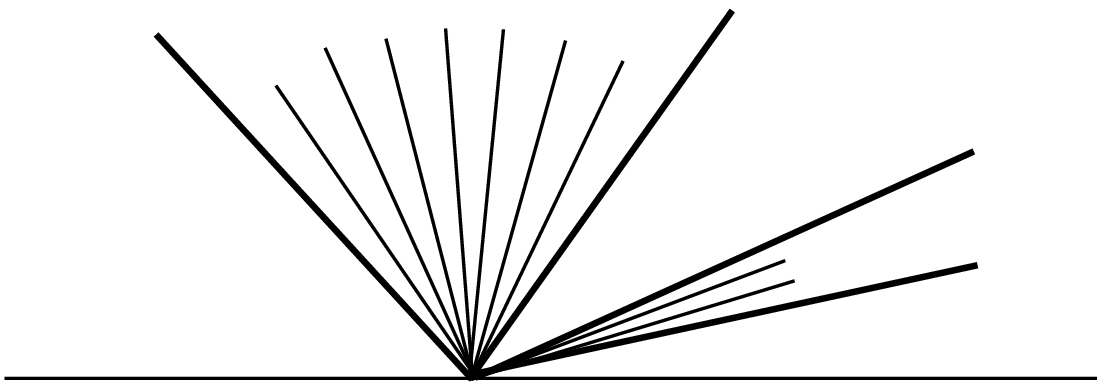}}
 \put(225,10){$(s_L,c_L)$}
 \put(230,-10){$(s_L,c_L)$}
 \put(282,-10){0}
  \put(310,-10){$(s_R,c_R)$}
  \put(330,5){$(s_R,c_R)$}
  \put(310,31){$(\overline{s},c_R)$}
  \put(235,55){$x=\sigma_1 t$}
   \put(300,55){$x=\sigma_2 t$}
   \put(350,35){$x=\sigma_3 t$}
\put(350,17){$x=\sigma_4 t$}
   \put(133,-10){$A$}
    \put(150,-10){$s_L$}
 \end{picture}
 \caption{Solution of Riemann problem (\ref{riemannpb}) with $s_L \geq s^*$ and
 $s_R < A$.}
 \label{rs2a}
 \end{figure}

\item {\bf Case 2b:} $s_R \geq A$\\

 Draw a line passing through the points $(s_R,f(s_R,c_R))$ and $(-\bar{a}_L(c_R),0)$. This line intersects the curve $f(s,c_L)$ at a point $\overline{s}$ (see Fig. \ref{rs2b}).
\begin{figure}[H]
 \begin{picture}(320,100)(-40,-10)
 \put(0,0){\includegraphics[height=3cm]{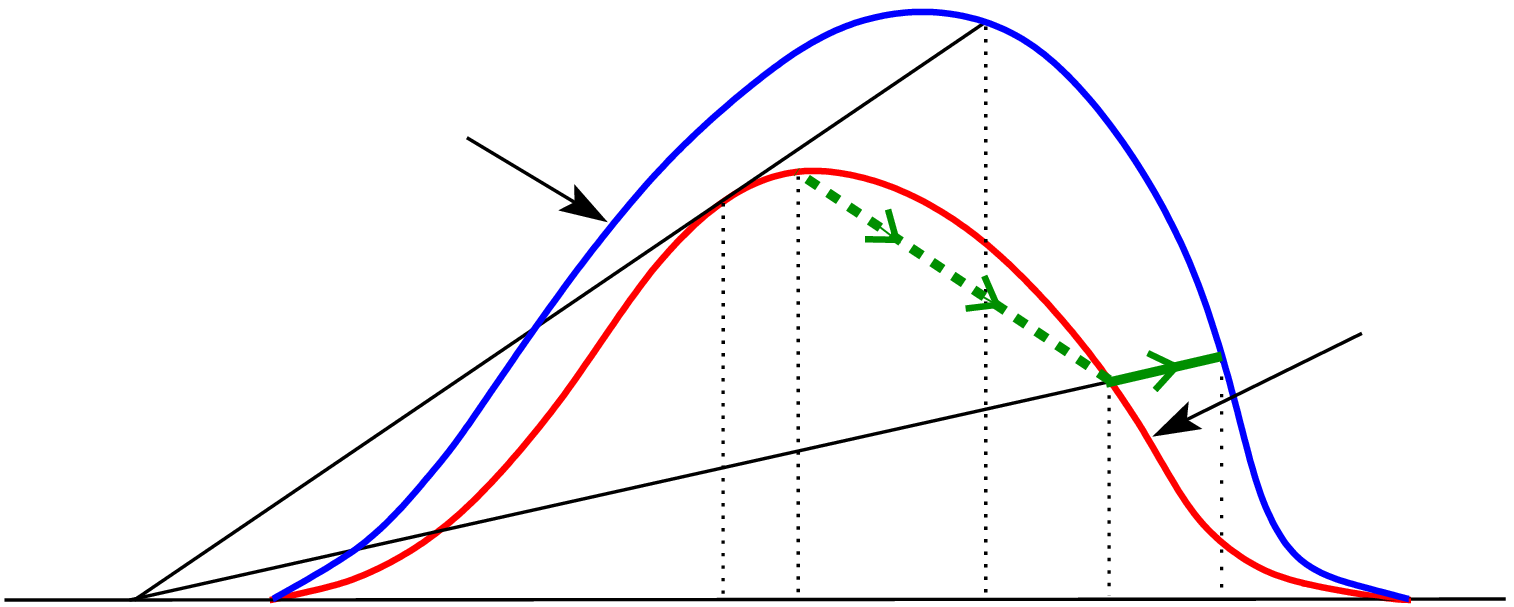}}
\put(40,72){$f(s,c_R)$}
 \put(190,45){$f(s,c_L)$}
 \put(108,-10){$s_L$}
\put(97,-10){$s^*$}
 \put(132,-10){$A$} 
 \put(157,-10){$\overline{s}$}
 \put(165,-10){$s_R$}
 \put(-7,-10){$-\bar{a}_L(c_R)$}
 \put(220,0){\includegraphics[height=1.8cm]{rs1bb.eps}}
 \put(225,15){$(s_L,c_L)$}
 \put(230,-10){$(s_L,c_L)$}
 \put(285,-10){0}
  \put(315,-10){$(s_R,c_R)$}
  \put(325,15){$(s_R,c_R)$}
  \put(283,28){$(\overline{s},c_L)$}
  \put(243,55){$x=\sigma_s t$}
   \put(340,55){$x=\sigma_c t$}
 \end{picture}
 \caption{Solution of Riemann problem (\ref{riemannpb}) with $s_L <  s^*$ and
 $s_R \geq A $.}
 \label{rs2b}
 \end{figure}
                                                                                
 (a)  Connect $(s_L,c_L)$ to $(\overline{s},c_L)$ by a $s$-wave along the curve $f(s,c_L)$,\\
 (b) Next connect  $(\overline{s},c_L)$ to $(s_R,c_R)$ by a $c$-wave with a speed
$$\sigma_c=\frac{f(s_R,c_R)}{s_R+\bar{a}_L(c_R)}=\frac{f(\overline{s},c_L)}{\overline{s}+\bar{a}_L(c_R)}. $$

    For example if $s_L < \overline{s}$ then the   solution is given by
\be 
(s(x,t),c(x,t)) = \left\{ \begin{array}{lll}
(s_L,c_L) &\mbox{if}& x<\sigma_s t, \\ (\overline{s},c_L) &\mbox{if}& \sigma_s
t < x < \sigma_c t, \\ (s_R,c_R) &\mbox{if}& x > \sigma_c t, \end{array} \right.
\ee
where
 $$\sigma_c=\frac{f(s_R,c_R)}{s_R+\bar{a}_L(c_R)}=\frac{f(\overline{s},c_L)}{\overline{s}+\bar{a}_L(c_R)}, \quad \sigma_s=\frac{f(\overline{s},c_L) - f(s_L,c_L)}{\overline{s}-s_L}.$$
\noi Note that $  \sigma_s < \sigma _c$ and $(s_L,c_L)$ is connected to
$(\overline{s},c_L)$ by a $s$-shock wave and $(\overline{s},c_L)$ is connected
to $(s_R,c_R)$ by a $c$-shock wave.
                                                                                
\end{itemize}

\section{Conservative finite volume schemes for the system of polymer flooding}
\label{finite difference}
    Let $h > 0$ and define the space grid points $x_{i+1/2}=ih,i \in \Z$ and
  for $\Delta t >0$ define the time discretization points  $t_n=n \Delta  t $ for all non-negative integer $n$. Let   
$\lambda=\frac{\Delta t}{h}$. A numerical scheme which is in conservative form for equation (\ref{polyeq1}) is given by

\begin{equation}
\begin{array}{l}
 (s_i^{n+1} - s_i^{n}) + \lambda ( F^n_{i+1/2} -  F^n_{i-1/2} )= 0,\\
(c_i^{n+1} s_i^{n+1} + a(c_i^{n+1})- c_i^n s_i^{n} - a(c_i^{n})) + 
\lambda (G^n_{i+1/2} -  G^n_{i-1/2})= 0
\label{finitevolumescheme}
\end{array}
 \end{equation}
 where the numerical flux $F^n_{i+1/2}$ and $ G^n_{i+1/2} $ are associated with the flux functions $f(s,c)$  and $g(s,c)=c f(s,c)$, and are functions of the left and right values of the saturation $s$ and the concentration $c$ at $x_{i+1/2}$:
 \[ F^n_{i+1/2} = F(s_i^n, c_i^n, s_{i+1}^n, c_{i+1}^n), \quad 
 G^n_{i+1/2} = G(s_i^n, c_i^n, s_{i+1}^n, c_{i+1}^n). \]
 The choice of the functions $F$ and $G$ determines the numerical scheme. We first present the new flux that we call DFLU, which is constructed as presented in the introduction. We compare it with the exact Riemann solver and show $L^\infty$ estimates for the associate scheme. Then we recall three other schemes to which to compare:  the upstream mobility flux and two centered schemes, Lax-Friedrichs's  and FORCE.
\subsection{The DFLU numerical flux}
\begin{figure}[htbp]   \begin{center}
\begin{picture}(350,70)(0,-5)
\thicklines
\put(70.,20){\parbox{100pt}{$\phi s_t + f(s, c_i^n)_x = 0$\\ $s(t_n) =
s_i^n$}}
\put(185.,20){\parbox{100pt}{$\phi s_t + f(s, c_{i+1}^n)_x = 0$\\ $s(t_n) =
s_{i+1}^n$}}
\put(0.,0.){\line(1,0){350}}
\put(0.,50.){\line(1,0){350}}
\put(175.,0.){\line(0,1){50}}
\put(60.,0.){\line(0,1){50}}
\put(300.,0.){\line(0,1){50}}
\put(160.,-10){$x_{i+1/2}$}
\put(45.,-10){$x_{i-1/2}$}
\put(285.,-10){$x_{i+3/2}$}
\put(0.,5){$t=t_n$}
\put(0.,55){$t=t_{n+1}$}
\end{picture}
\end{center}
\label{fluxdiscr}
\caption{The flux functions $f(\cdot,c)$ is discontinuous in c at the
discretization points.}
\end{figure}
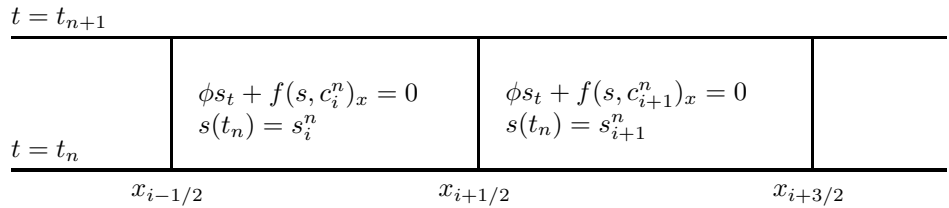

The DFLU flux is an extension of the Godunov scheme that we proposed and analyze in \cite{AdiJafGow04} for scalar conservations laws with a flux function discontinuous in space.
As the second eigenvalue $\lambda^c$ of the system is always non-negative we define
\begin{equation}
G^n_{i+1/2}=c^n_i\,F^n_{i+1/2}. 
 \label{DFLUg} 
\end{equation}. 
Now the choice of the numerical scheme depends on the choice
of $F^n_{i+1/2}$. To do so we treat $c(x,t)$
in $f(s,c)$  as  a known function which may be discontinuous at the space discretization points.
Therefore on the border of each rectangle $(x_{i-1/2},x_{i+1/2}) \times (t_n,t_{n+1})$, we
 consider the conservation law:
\begin{equation}
s_t+f(s,c_i^n)_x=0 
\end{equation}
 with initial condition $s(x,0)=s_i^0$ for $x_{i-1/2}< x< x_{i+1/2}$.
(see Fig.\ref{fluxdiscr}).

Extending the idea of \cite{AdiJafGow04},we define the DFLU flux  as
\be
\begin{array}{lcll}
 F_{i+1/2}^{n} & = & F^{DFLU}(s_i^n, c_i^n, s_{i+1}^n, c_{i+1}^n) \\
             & = & \min\{ f(\min\{s_i^n,\theta_i^n\},c_i^n),f(\max\{s_{i+1}^n,\theta_{i+1}^n\},c^n_{i+1})\},
\end{array} 
\label{DFLU}
\ee
where
$ \theta_i^n =\mbox{ argmax } f(\cdot,c_i^n) $.\\
                                                                         
\noindent {\bf Remarks:} \\
1) Suppose $c_i^n=c_0$, a constant for all $i$,then it is easy to see that $c_i^{n+1}=c_0$ for all $i$.\\
2) Suppose $s \rightarrow f(s,c)$ is an increasing function (case without gravity)  then $\theta_i^n=1$ for all $i$ and from (\ref{DFLU}) we have $ F^n_{i+1/2} = f(s_i^n,c_i^n)$  and the finite difference scheme (\ref{finitevolumescheme}) becomes
\be
\begin{array}{rll}
 s_i^{n+1} &= & s_i^n - \lambda (f(s_i^n,c_i^n) -f(s_{i-1}^n,c_{i-1}^n ))\\
c_i^{n+1} s_i^{n+1} + a(c_i^{n+1})&=& c_i^n s_i^{n} + a(c_i^{n}) -
 \lambda (c_i^n\,f(s_i^n,c_i^n) -c_{i-1}^n\,f(s_{i-1}^n,c_{i-1}^n ))
\label{finitediffeqn2}
\end{array}
\ee
 which  is nothing but the standard upwind scheme.

\subsection{Comparison of the  DFLU flux with the flux given by an exact Riemann solver}
\label{comparison}

     Now we would like to compare the exact Godunov flux $F_{i+1/2}^G $ with our DFLU flux $F^{DFLU}_{i+1/2}$ defined by (\ref{DFLU}). For sake of brevity we considered only the case $c^n_i \geq c^n_{i+1}$. The opposite case can be considered similarly. We discuss the cases considered in section \ref{Riemann}.

\noi {\bf Case 1a:} $ s_i < s^*, s_{i+1} <  B$. See  Fig. \ref{figR1a}. 
   In this case $F_{i+1/2}^G=f(s_i,c_i)=F^{DFLU}_{i+1/2}$.\\

\noi {\bf Case 1b:} $ s_i < s^*, s_{i+1} \geq B$. See  Fig. \ref{rs1b}. \\
   Then 
$F_{i+1/2}^G= \left\{ \begin{array}{lll}
   f(\overline{s},c_i) &\mbox{if}& \sigma_s < 0 \\ f(s_i,c_i) &\mbox{if}& \sigma_s  \geq 0 
\end{array} \right. $
where $ \sigma_s=\dfrac{f(\overline{s},c_i) - f(s_i,c_i)}{\overline{s}-s_i}$. On the other hand the DFLU flux gives   $ F^{DFLU}_{i+1/2}=\min \{f(s_i,c_i),f(\max \{s_{i+1},\theta_{i+1} \},c_{i+1}) \} $.
Therefore in this case the Godunov flux may not be  same as the DFLU flux.\\

\noi {\bf Case 2a:} $s_i \geq s^*,s_{i+1} \leq A $. See Fig.\ref{rs2a}. Then
\[
F_{i+1/2}^G = \left\{ \begin{array}{lll}
   f(\theta_i,c_i) &\mbox{if} &  s_i > \theta_i   \\ f(s_i,c_i) &\mbox{if}& s_i \leq \theta_i
\end{array} \right.
 =f(\min \{s_i,\theta_i \},c_i) =F_{i+1/2}^{DFLU}.
\]
\noi {\bf Case 2b:}$s_i \geq s^*,s_{i+1} > A $. See Fig.\ref{rs2b}. \\
Then $F_{i+1/2}^G= \left\{ \begin{array}{lll}
   f(\overline{s},c_i) &\mbox{if}&  \sigma_s < 0   \\ f(s_i,c_i) &\mbox{if}& \sigma_s \geq 0
\end{array} \right. $
where $ \sigma_s=\dfrac{f(\overline{s},c_i) - f(s_i,c_i)}{\overline{s}-s_i}$. \\
The DFLU flux is $F^{DFLU}_{i+1/2}= \min\{ f(\min\{s_i,\theta_i\},c_i),f(\max\{s_{i+1},\theta_{i+1}\},c_{i+1})\}$.
In this case these two fluxes are not equal, for example when $\sigma_s <0$.

One can actually observe that the Godunov flux can actually be calculated with  the following compact formula:

\noi {\bf Case 1:  $s_i < s^*_i$.}
\[
F_{i+1/2}^G = \!\! \left\{ \begin{array}{lllll}
  \!\! f(s_i,c_i) &\! \mbox{if}\, f_s(s_{i+1},c_{i+1})\geq 0 
  \mbox{ or }\, \dfrac{f(s_{i+1},c_{i+1})}{s_{i+1}+\bar{a}_L(c_{i+1})} \geq \dfrac{f(s_i,c_i)}{s_i+\bar{a}_L(c_{i+1})},\\
   \!\! \min(f(s_i,c_i),f(\overline{s}_i,c_i))&  \!\mbox{otherwise},
\end{array} \right.
\]
where $\overline{s}_i$ is given by
$\dfrac{f(s_{i+1},c_{i+1})}{s_{i+1}+\bar{a}_L(c_{i+1})}=\dfrac{f(\overline{s}_i,c_i)}{\overline{s_i}+\bar{a}_L(c_{i+1})}$.

\noi {\bf Case 2:  $s_i \ge s^*_i$.}
\[
F_{i+1/2}^G = \left\{ \begin{array}{llll}
   \!\! f(\min(s_i,\theta_i),c_i) &\! \mbox{if} \, f_s(s_{i+1},c_{i+1})\geq 0 
   \mbox{ or }\, \dfrac{f(s_{i+1},c_{i+1})}{s_{i+1}+\bar{a}_L(c_{i+1})} \geq \dfrac{f(s^*_i,c_i)}{s^*_i+\bar{a}_L(c_{i+1})},&\\
  \!\!  \min(f(s_i,c_i),f(\overline{s}_i,c_i)) & \! \mbox{otherwise},
\end{array} \right.
\]
where $\overline{s}_i$ is given by
$\dfrac{f(s_{i+1},c_{i+1})}{s_{i+1}+\bar{a}_L(c_{i+1})}=\dfrac{f(\overline{s}_i,c_i)}{\overline{s}_i+\bar{a}_L(c_{i+1})}$.

\subsection{$L^{\infty}$ and TVD bounds for the DFLU scheme}
We show first $L^{\infty}$ bounds, and TVD bounds will follow immediately. Let 
$\displaystyle{M=\sup_{s,c}\{f_s(s,c),\frac{f(s,c)}{s+a'(c)}\}}$. 
\begin{lemma}
Let $s_0$ and $c_0 \in  L^\infty (\R, [0,1])$ be the initial data and let
$\{s_i^n\}$ and $\{c_i^n\}$ be the corresponding solution calculated
by the finite volume scheme (\ref{finitevolumescheme}) using the DFLU flux (\ref{DFLUg}),  (\ref{DFLU}).
When  $\la M \leq 1$  then
\be
\begin{array}{l}
0 \leq s_i^n \leq 1 \; \mbox{ for all } i, n,\\
||c^n||_{\infty}  \leq    ||c^{n-1}||_{\infty} \mbox{ where } ||c^n||_{\infty}=\sup_i|c_i^n|.
\end{array}
\label{scinfinity}
\ee
\end{lemma}
{\bf Proof:} Since $0 \leq s_0 \leq 1$ and hence for all $i, \; 0 \leq s_i^0
\leq 1.$ By induction, assume that (\ref{scinfinity}) holds for all $n$. Let 
$$\begin{array}{rll}
 s_i^{n+1} &=& s_i^{n} - \lambda ( F^n_{i+1/2} -  F^n_{i-1/2} )\\
           & =& H(s_{i-1}^n,s_i^n,s_{i+1}^n,c_{i-1}^n,c_i^n,c_{i+1}^n)
\end{array}
$$
By (\ref{DFLU}),it is easy to check that if $\lambda M \leq 1$, then
$H=H(s_1,s_2,s_3,c_1,c_2,c_3)$ is an increasing function in $s_1,s_2,s_3$ and by the hypothesis on $f$,$H(0,0,0,c_1,c_2,c_3)=0, H(1,1,1,c_1,c_2,c_3)=1$. Therfore

$$\begin{array}{rll}
0 &=& H(0,0,0,c_{i-1}^n,c_i^n,c_{i+1}^n) \\
&\leq & H(s_{i-1}^n,s_i^n,s_{i+1}^n,c_{i-1}^n,c_i^n,c_{i+1}^n)=s_i^{n+1}\\
& \leq &  H(1,1,1,c_{i-1}^n,c_i^n,c_{i+1}^n)=1.
\end{array}
$$
This proves $0  \leq s_i^{n+1}  \leq 1$.

     To prove the boundness of $c$, consider   
$$ (c_i^{n+1} s_i^{n+1} + a(c_i^{n+1})- c_i^n s_i^{n} - a(c_i^{n})) +
\lambda (G^n_{i+1/2} -  G^n_{i-1/2})= 0.$$
By adding and subtracting the term $c_i^ns_i^{n+1}$ to the second equation
of (\ref{finitevolumescheme}) and by  substituting  first equation we can rewrite the second equation as
$$ c_i^{n+1}(s_i^{n+1}+a'(\xi_i^{n+1/2}))-c_i^n (s_i^{n+1}+a'(\xi_i^{n+1/2}))+ \lambda  F^n_{i-1/2} (c_i^n-c_{i-1}^n)=0$$
where $ a(c_i^{n+1})- a(c_i^{n})=a'(\xi_i^{n+1/2}) (c_i^{n+1}-c_i^n)$ for some
$\xi_i^{n+1}$ between $c_i^{n+1}$ and $c_i^n$.
   This is equivalent to 
$$ c_i^{n+1}=c_i^n - \lambda \frac{ F^n_{i-1/2}}{(s_i^{n+1}+a'(\xi_i^{n+1/2}))} (c_i^n-c_{i-1}^n)$$
which is the scheme written in the non-conservative form.
Let $b_i^n=\lambda \dfrac{ F^n_{i-1/2}}{(s_i^{n+1}+a'(\xi_i^{n+1/2}))}$ then
\[
 c_i^{n+1} = (1-b_i^n)c_i^n+b_i^nc_{i-1}^n
          \leq \max\{c_i^n,c_{i-1}^n\} \mbox{ if } b_i^n \leq 1.
\]
This proves the second inequality.  \cqfd\\

Since $c_i^{n+1}$ is a convex combination
of $c_i^n$ and $c_{i-1}^n$ if $\lambda M \leq 1$, then we obtain the following total variation diminishing property for $c_i^{n}$:
\begin{lemma}
 Let $\{c_i^n\}$ be the solution  calculated
by the finite volume scheme (\ref{finitevolumescheme}), (\ref{DFLUg}), (\ref{DFLU}).
When  $\la M \leq 1$  then
\[
\sum_i{|c_i^{n+1} - c_{i-1}^{n+1}|}  \leq  \sum_i|c_i^{n} -c_{i-1}^{n}|\,\, \mbox{ for all n}.
\]
\end{lemma}
Note that the saturation itself is not TVD because of the discontinuity of $f$, and that the above proof applies also to the usptream mobility flux presented below.

\subsection{The upstream mobility flux}

 Petroleum engineers have designed, from physical
considerations, another numerical flux
called the upstream mobility flux. It is an ad-hoc flux for two-phase
flow in porous media which
corresponds to an approximate solution to the Riemann problem. For this flux $G^n_{i+1/2}$ is given again by (\ref{DFLUg}) and $F_{i+\frac{1}{2}}^n$ is given by
\[
\begin{array}{l}
F_{i+\frac{1}{2}}^n=F^{UM}(s_i^n,c_i^n,s_{i+1}^n,c_{i+1}^n)= \displaystyle{\frac{1}{\phi} \,
         \frac{\lambda_1^*}{\lambda_1^* + \lambda_2^*}
         [ q + (c_1-c_2)\lambda_2^* ]} ,\\
\lambda^*_\ell=
\left\{ \begin{array}{ll} \lambda_\ell(s^n_i,c_i^n) & \mbox{if }
         q+(g_\ell-g_i)\lambda_\ell^* >0, \;i=1,2, i\neq \ell,\\[3mm]
        \lambda_\ell(s_{i+1}^n,c_{i+1}^n) & \mbox{if }
        q+(g_\ell-g_i)\lambda_\ell^* \leq0, \; i=1,2, i\neq \ell,
\end{array}
\right. \ell=1,2.
\end{array}
\]                                                               

\subsection{The Lax-Friedrichs flux}
In this case fluxes are given by
\[
\begin{array}{llll}
F^n_{i+1/2}&=& \frac{1}{2} [f(s^n_{i+1},c^n_{i+1})+f(s^n_i,c^n_i)-\dfrac{(s^n_{i+1}-s^n_i)}{\lambda} ]\\
G^n_{i+1/2}&=&\frac{1}{2} [c^n_{i+1}f(s^n_{i+1},c^n_{i+1})+c^n_i
f(s^n_i,c^n_i)-\dfrac {(c^n_{i+1} s^n_{i+1}+a(c^n_{i+1})-c^n_is^n_i-a(c^n_i))}
{\lambda} ]
\end{array}
\]

\subsection{The FORCE  flux}
This flux \cite{Toro99,AdiGowJaf09}, introduced by E. F. Toro, is an average of the Lax-Friedrichs and Lax-Wendroff flux. It is defined by
\[
\begin{array}{llll}
F^n_{i+1/2}&=& \frac{1}{4} [f(s^n_{i+1},c^n_{i+1})+f(s^n_i,c^n_i)+2 f(s^{n+1/2}_i)-\dfrac{(s^n_{i+1}-s^n_i)}{\lambda} ]\\[0.3cm]
G^n_{i+1/2}&=&\frac{1}{4} [c^n_{i+1}f(s^n_{i+1},c^n_{i+1})+c^n_i
f(s^n_i,c^n_i)+2 c^{n+1/2}_i f(s^{n+1/2}_i,c^{n+1/2}_i)\\[0.2cm]
&& -\dfrac {(c^n_{i+1} s^n_{i+1}+a(c^n_{i+1})-c^n_i s^n_i-a(c^n_i))}{\lambda}]\\
\end{array}
\]
 where
$$s^{n+1/2}_i=\frac{(s^n_{i+1}+s^n_i)}{2}-\frac {\lambda}{2}(f(s^n_{i+1},c^n_{i+1})-f(s^n_i,c^n_i))$$
and
\[
\begin{array}{llll}
 s^{n+1/2}_ic^{n+1/2}_i+a(c^{n+1/2}_i)&=&\dfrac{(s^n_{i+1}c^n_{i+1}+s^n_ic^n_i)}{2}+\dfrac{(a(c_{i+1}^n)+a(c_i^n))}{2}\\[0.3cm]
&&\qquad \qquad -\frac{\lambda}{2}(c^n_{i+1}f(s^n_{i+1},c^n_{i+1})-c^n_if(s^n_i,c^n_i)).
\end{array}
\]

\section{Numerical experiments}
\label{numericalresults}
  To evaluate the performance of the DFLU scheme we first compare its results to an exact solution and evaluate convergence rates, and then  compare it with other standard numerical
 schemes already mentioned in the previous section, that are the Godunov, upstream mobility, Lax-Friedrichs and FORCE schemes.  
 \subsection{Comparison with an exact solution}
In this section we compare the calculated and exact solutions of two Riemann problems. We consider the following functions 
\be  f(s,c)=s(4-s)/(1+c),  \quad a(c)=c. \label{flux1} \ee
Note that $f(0,c)=f(4,c)=0$ for all $c$ and the interval for $s$ is $[0,4]$ instead of $[0,1]$. This choice of $f$, which does not correspond to any physical reality, was done in order to try to have a large difference between the Godunov and the DFLU flux (see second experiment below).

In a first experiment the initial condition is
\be s(x,0) = \left\{ \begin{array}{lll}
2.5 &\mbox{if}& x<.5, \\ 1 &\mbox{if}& x>.5 \end{array} \right.  , \quad
c(x,0) = \left\{ \begin{array}{lll}
.5 &\mbox{if}& x<.5, \\ 0 &\mbox{if}& x>.5. \end{array} \right.  
\label{IC1} \ee
  These $f$ and initial data correspond to 
   the case 2a in sections \ref{Riemann} and \ref{comparison} where the DFLU flux coincides with the Godunov flux: $F^{DFLU}(s_L,s_R,c_L,c_R)=F^{G}(s_L,s_R,c_L,c_R)$ with 
   $s^*=1.236, A=2.587,\overline{s}=.394$.  The exact  solution of the Riemann problem at a time $t$ is given by
\be s(x,t) = \left\{ \begin{array}{lll}
2.5 &\! \mbox{if}&\! x<.5+\sigma_1 \,t, \\ \frac{1}{2}(4-1.5(\frac{x-.5}{t}))&\!\mbox{if}&\! .5+\sigma_1\,t <x < .5 +\sigma_c
\, t \\  \overline{s}=.394 &\!\mbox{if}&\! .5+\sigma_c\,t < x <.5+\sigma_2 \,t \\
1. & \!\mbox{if} &\! x > \sigma_2 t+.5
  \end{array} \right.  , \;\;
c(x,0) = \left\{ \begin{array}{lll}
.5 & \!\!\! \mbox{if}&  \!\!\!  x<.5+\sigma_c t, \\ 0. &  \!\!\! \mbox{if}&  \!\!\!  x>.5+ \sigma_c t. \end{array} \right.  
\ee
where $\sigma_1=f_s(s_L,c_L)=-2/3$,
 $\sigma_c=f_s(s^*,c_L)=\dfrac{f(s^*,c_L)}{s^*+\bar{a}_L(c_R)}=\dfrac{f(\overline{s},c_R)}{\overline{s}+\bar{a}_L(c_R)}=1.018$ and $ \sigma_2=\dfrac{f(\overline{s},c_R) - f(s_R,c_R)}{\overline{s}-s_R}=2.606$.

\begin{figure}[htbp]
 \includegraphics[width=7.2cm]{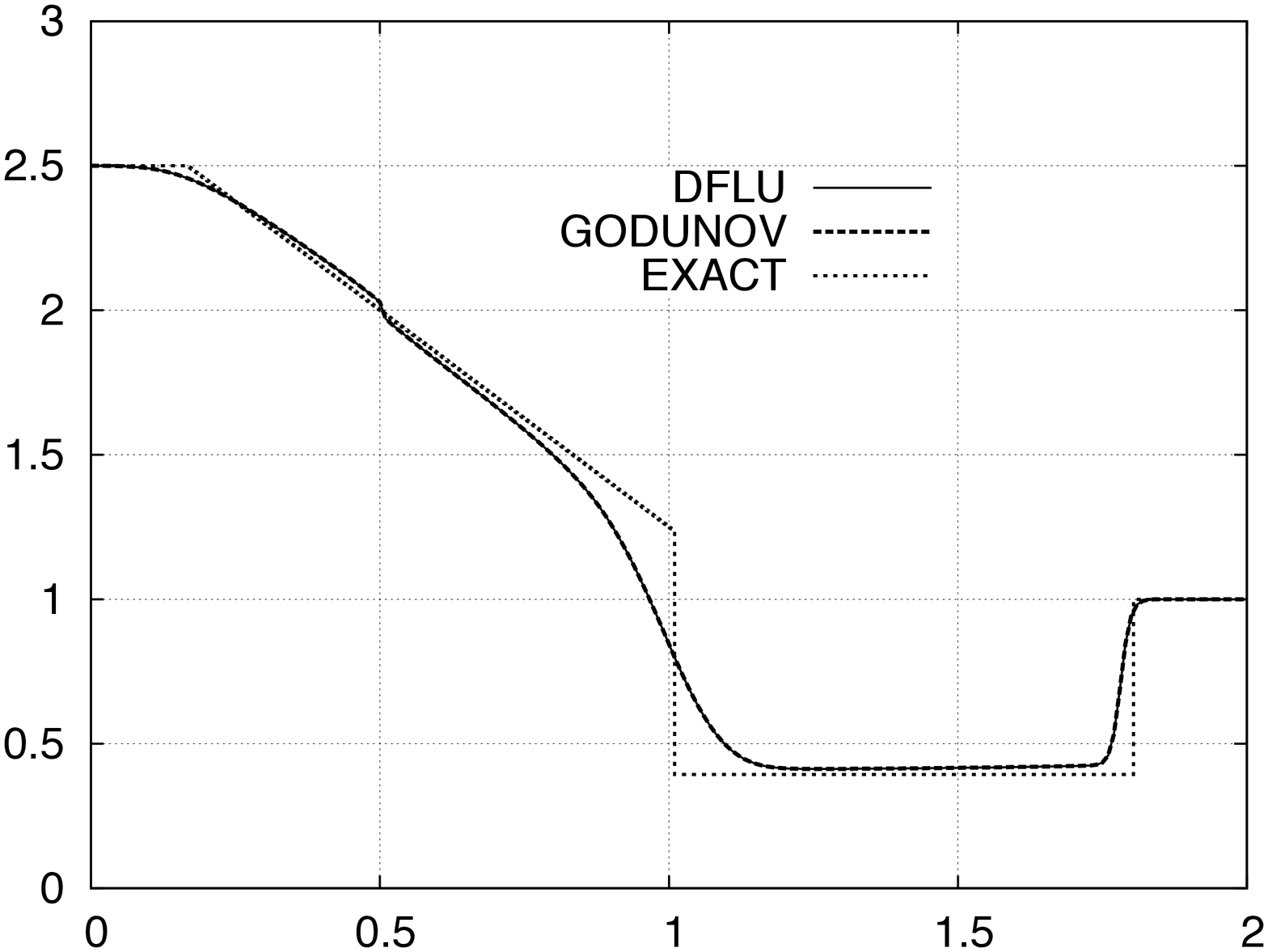} \hspace*{0.3cm}
 \includegraphics[width=7.2cm]{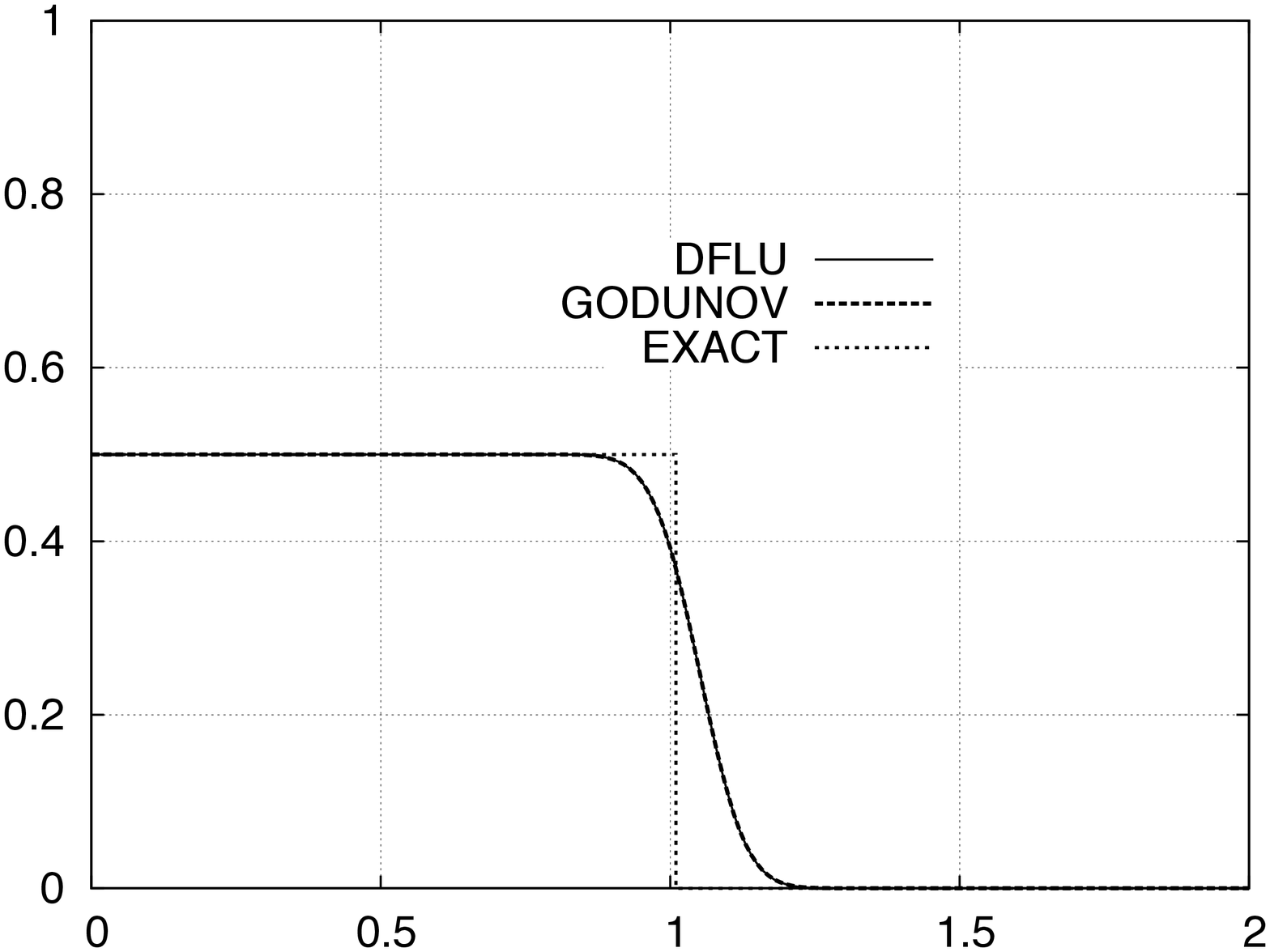}
\caption{Comparison with exact solution of Riemann problem (\ref{flux1}), (\ref{IC1}):  $s$ (left) and $c$ (right) at $t=.5$ for $h=1/100, \lambda=1/4$.} 
\label{fig2aexp1}
\end{figure}

\begin{figure}[H]
 \includegraphics[width=7.2cm]{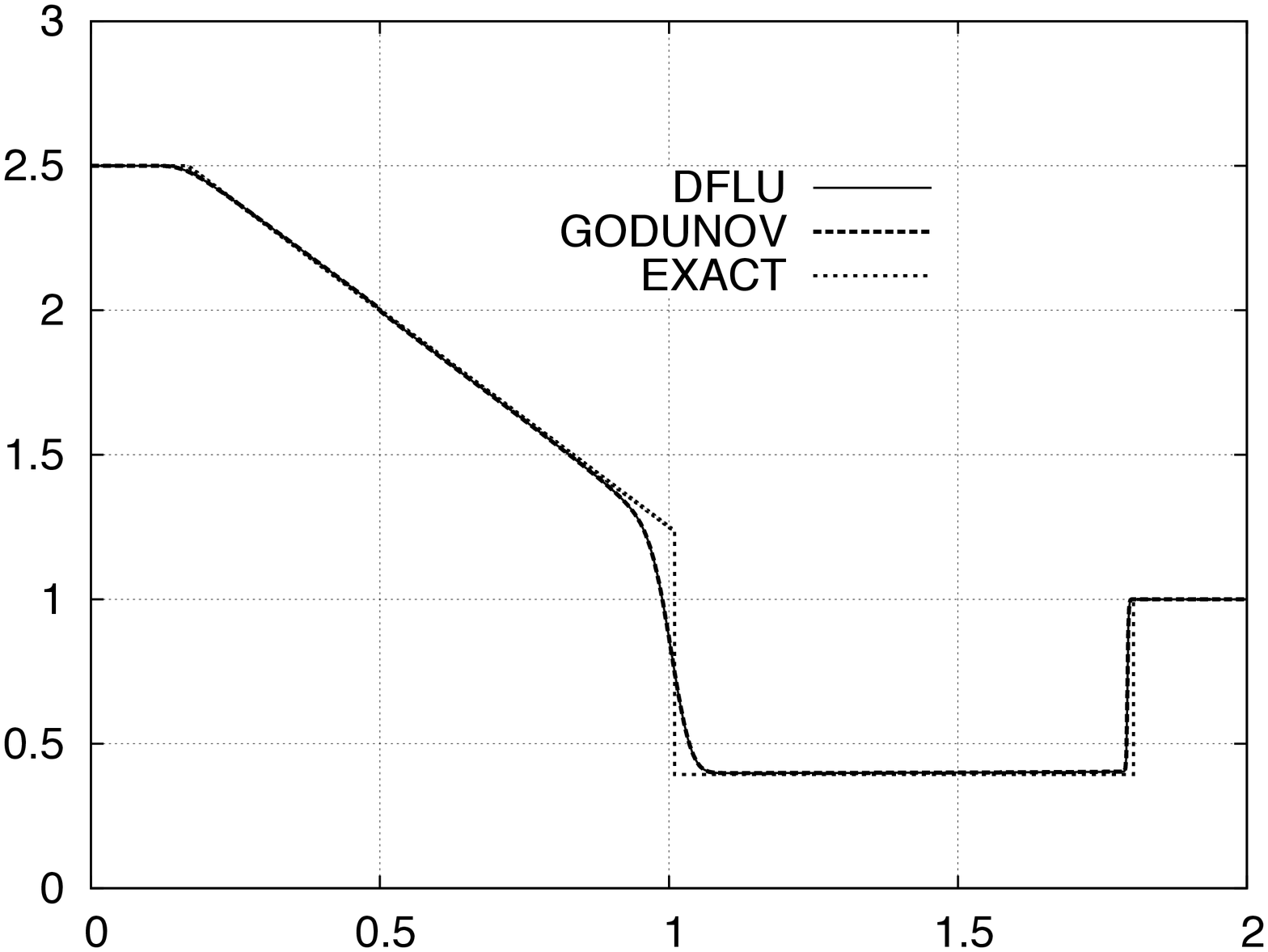} \hspace*{0.3cm}
 \includegraphics[width=7.2cm]{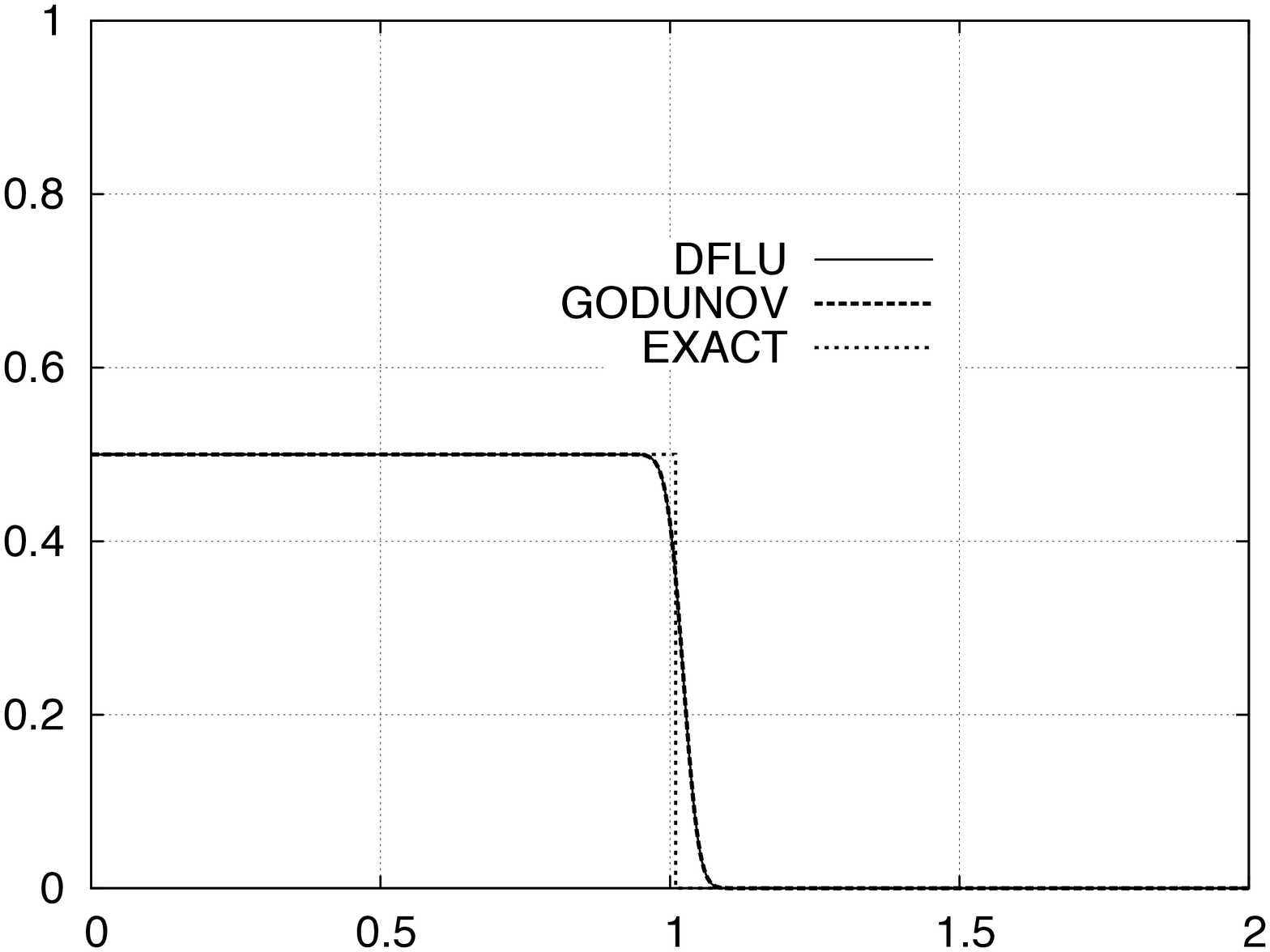}
\caption{Comparison with exact solution of Riemann problem (\ref{flux1}), (\ref{IC1}):  $s$ (left) and $c$ (right) at $t=.5$ for $h=1/800, \lambda=1/4$.} 
\label{fig2aexp11}
\end{figure}

 Figs. \ref{fig2aexp1} and \ref{fig2aexp11} verify that the DFLU and Godunov schemes give coinciding results.  
 As expected both  schemes are diffusive at $c$-shocks as well as at $s$-shocks but
 as the mesh size goes to zero calculated solutions are getting closer  
 to the exact solution (see Fig.\ref{fig2aexp11}).  Table \ref{table1} shows $L_1$ errors for $s$ and  $c$ and the convergence rate $\alpha$. Calculations are done with $\lambda=\frac{1}{4} (M=4)$, that is the largest time step allowed by the CFL condition.
\begin{table}[hbtp]
\begin{center}
\begin{tabular}{|l||c|c||c|c|}\hline
$h$& Godunov,$||s-s_h||_{L^1}$&$\alpha$&DFLU,$||s-s_h||_{L^1}$& $\alpha$ \\
\hline
1/50&.2373&&.2372&\\
\hline
1/100&0.15134&0.6489&0.1506&0.655\\
\hline
1/200& 9.6868 $\times 10^{-2}$&0.6437&9.6868 $\times 10^{-2}$&0.6366  \\
\hline
1/400&6.4228 $\times 10^{-2}$&0.5928&6.4228 $\times 10^{-2}$ &  0.5928\\
\hline
1/800&4.2198 $\times 10^{-2}$&0.606&4.2197 $\times 10^{-2}$ &0.606  \\
\hline
\end{tabular}\\[0.3cm]
\begin{tabular}{|l||c|c||c|c|}\hline
$h$& Godunov,$||c-c_h||_{L^1}$&$\alpha$&DFLU,$||c-c_h||_{L^1}$& $\alpha$\\
\hline
1/50& 6.3796 $\times 10^{-2}$&&6.3796 $\times 10^{-2}$&\\
\hline
1/100&4.1630 $\times 10^{-2}$&0.6158& 4.1630 $\times 10^{-2}$&0.6158 \\
\hline
1/200& 2.6669 $\times 10^{-2}$&0.6424&2.6669 $\times 10^{-2}$ &0.6424 \\
\hline
1/400&1.7398 $\times 10^{-2}$&0.6162&1.7398 $\times 10^{-2}$ &0.6162  \\
\hline
1/800&1.1522 $\times 10^{-2}$&0.5945&1.1522 $\times 10^{-2}$&0.5945  \\
\hline
\end{tabular}
\end{center}
\caption{Riemann problem (\ref{flux1}), (\ref{IC1}): $L^1$-errors between exact and calculated solutions at $t=.5$}
\label{table1}
\end{table}

Now we want to have an experiment where the DFU flux differs from the Godunov flux. Therefore we now consider the Riemann problem with initial data
\be s(x,0) = \left\{ \begin{array}{lll}
2.3 &\mbox{if}& x<.5, \\ 3.2 &\mbox{if}& x>.5, \end{array} \right.  , \quad
c(x,0) = \left\{ \begin{array}{lll}
.5 &\mbox{if}& x<.5, \\ 0 &\mbox{if}& x>.5. \end{array} \right.  
\label{IC2}
\ee
This initial data corresponds to case 2b of  sections \ref{Riemann} and \ref{comparison} with
 $c_R=0$, $s^*=1.236$.
     In this case, the  exact solution of  the Riemann problem at a time $t$ is given by
\[ s(x,t) = \left\{ \begin{array}{lll}
s_L=2.3 &\mbox{if}& x<.5+\sigma_s \,t \\
 \overline{s}=2.7536&\mbox{if}& .5+\sigma_s\,t <x < .5 +\sigma_c t,\\
s_R=3.2 &\mbox{if} & x > \sigma_c t+.5
  \end{array} \right.  , \quad
c(x,0) = \left\{ \begin{array}{lll}
.5 &\mbox{if}& x<.5+\sigma_c t, \\ 0. &\mbox{if}& x>.5+ \sigma_c t, \end{array} \right.  
\]
where $\sigma_s=\dfrac{f(s_L,c_L)-f(\overline{s},c_L)}{s_L-\overline{s}}=-.702$,
 and  $\sigma_c=\dfrac{f(s_R,c_R)}{s_R+\bar{a}_L(c_R)}=0.609$.       
 
\begin{figure}[H]
 \includegraphics[width=7.2cm]{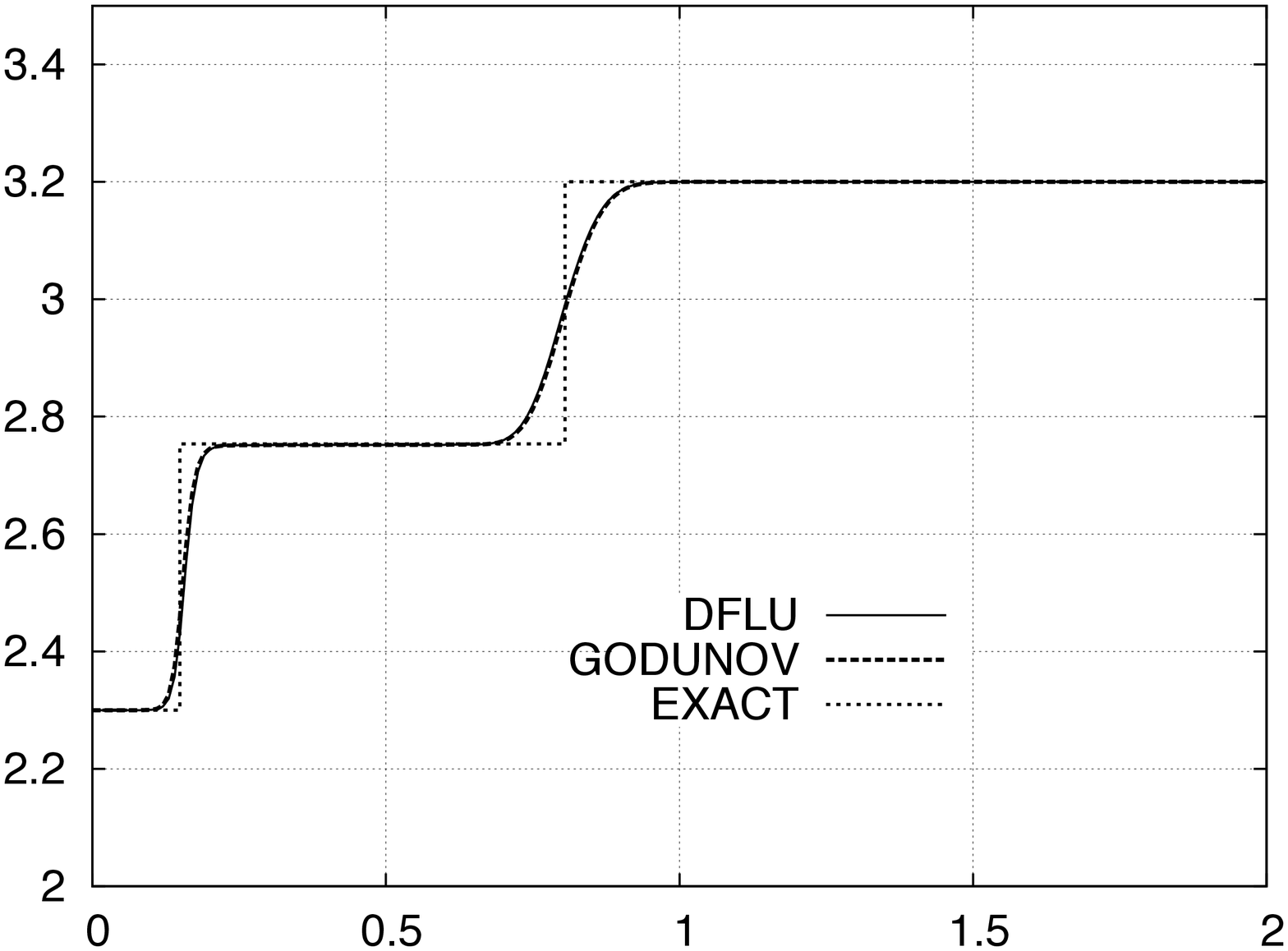}
 \includegraphics[width=7.2cm]{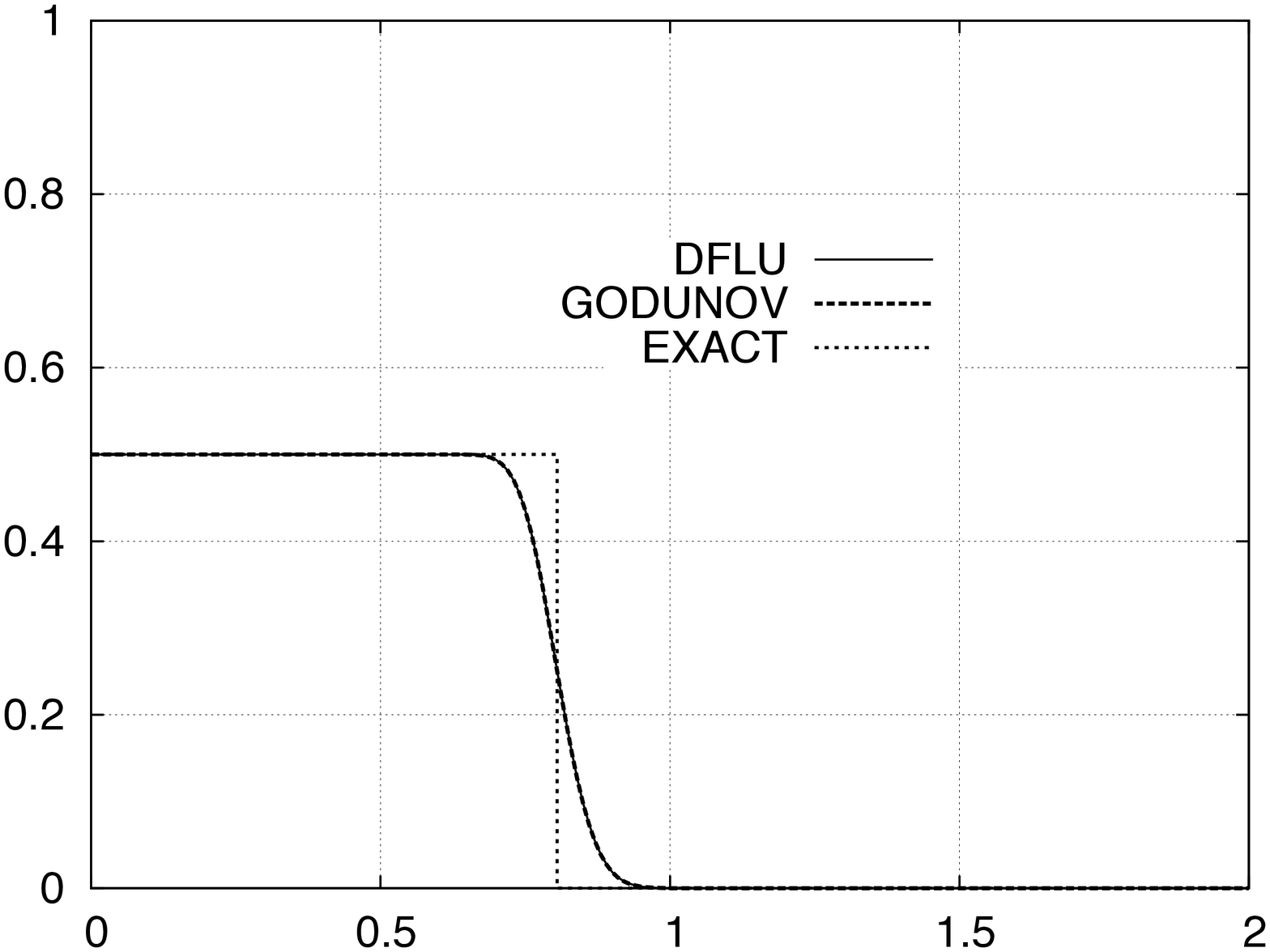}
\caption{Comparison with exact solution of Riemann problem (\ref{flux1}), (\ref{IC2}):  $s$ (left) and $c$ (right) at $t=.5$ for $h=1/100, \lambda=1/4$.} 
\label{fig2bexp2}
\end{figure}

\begin{figure}[H]
 \includegraphics[width=7.2cm]{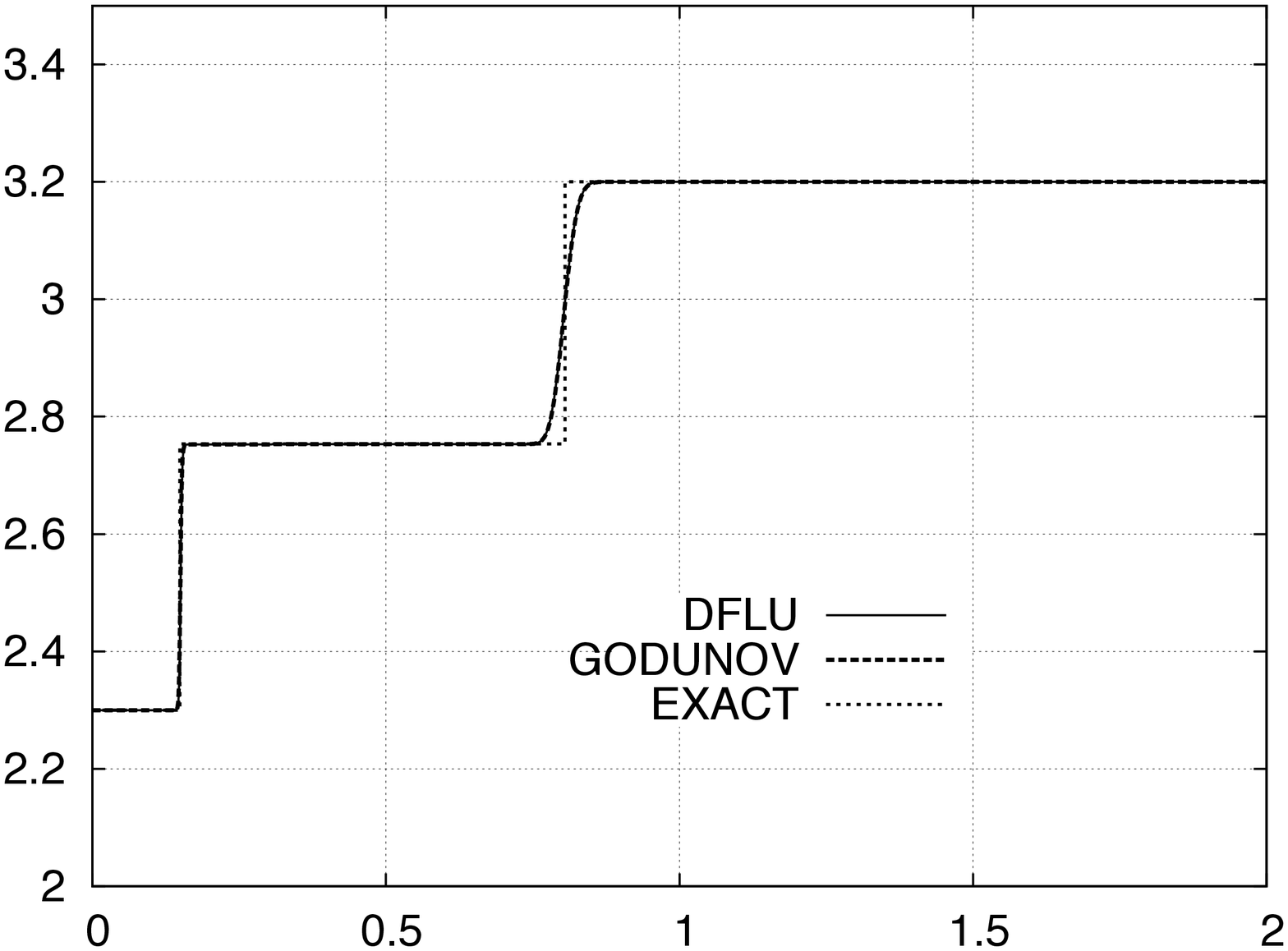}
 \includegraphics[width=7.2cm]{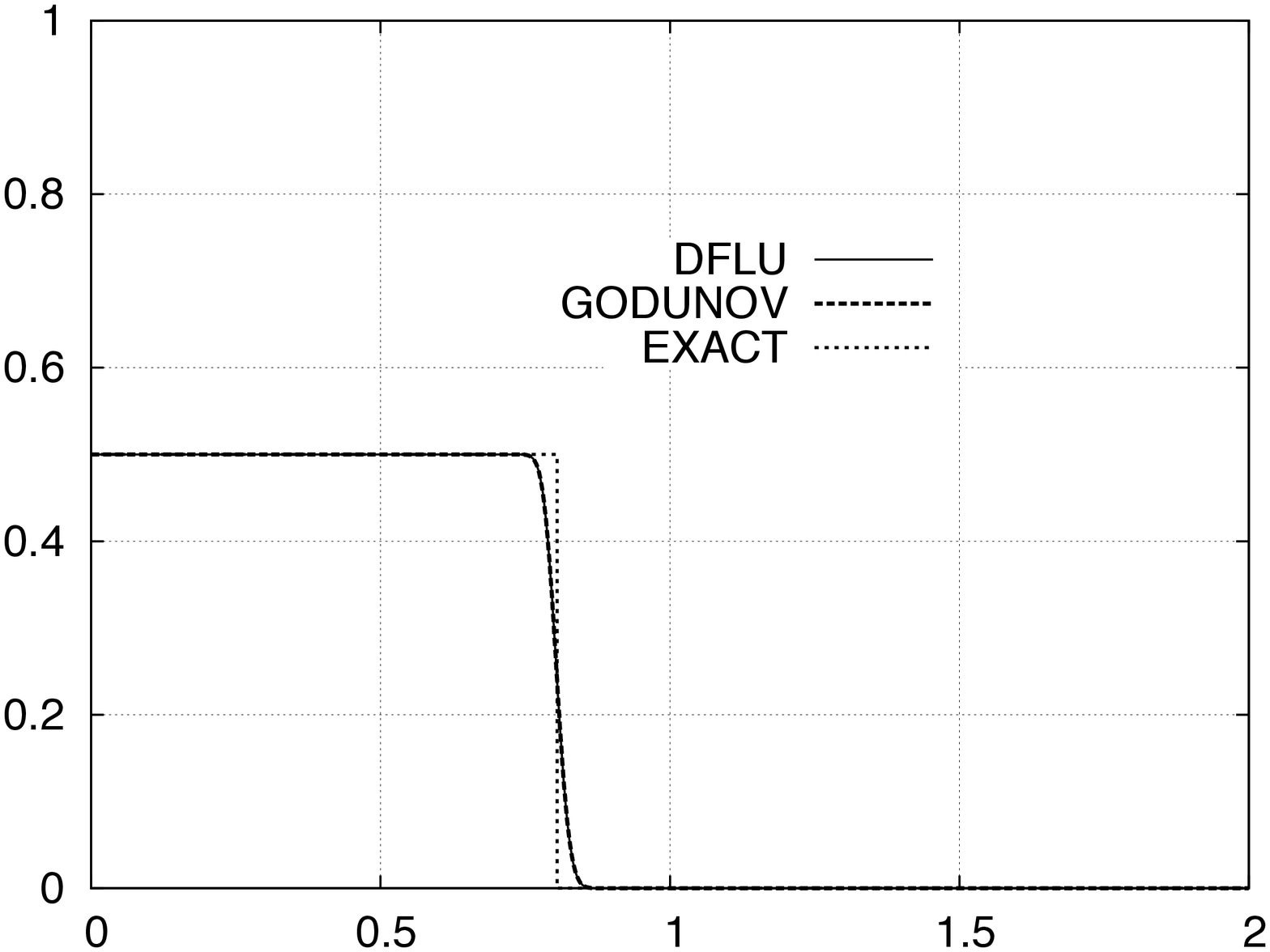}
\caption{ Comparison with exact solution of Riemann problem (\ref{flux1}), (\ref{IC2}):  $s$ (left) and $c$ (right) at $t=.5$ for $h=1/800, \lambda=1/4$.} 
\label{fig2bexp21}
\end{figure}

Figs. \ref{fig2bexp2} and \ref{fig2bexp21} show the comparison of the results obtained with the DFU and Godunov fluxes with the exact solution. The solution obtained with the DFU and Godunov flux are very close even if they do not coincide actually. Table \ref{table2} shows $L_1$ errors for $s$ and  $c$ and the convergence rate $\alpha$. Calculations are done with $\lambda=\frac{1}{4} (M=4)$, that is the largest time step allowed by the CFL condition.

\begin{table}[H]
\begin{center}
\begin{tabular}{|l||c|c||c|c|}\hline
$h$& Godunov,$||s-s_h||_{L^1}$&$\alpha$&DFLU,$||s-s_h||_{L^1}$&$\alpha$\\
\hline
1/50&0.10246&&0.10373&\\
\hline
1/100&5.7861 $\times 10^{-2}$ &0.8243 &5.8731 $\times 10^{-2}$&0.8206\\
\hline
1/200& 3.2849 $\times 10^{-2}$&0.81674&3.3259 $\times 10^{-2}$& 0.8203 \\
\hline
1/400&1.9152 $\times 10^{-2}$&0.7785&1.9353 $\times 10^{-2}$ & 0.7811 \\
\hline
1/800&1.1489 $\times 10^{-2}$&0.7370&1.1571 $\times 10^{-2}$ &  0.7420\\
\hline
\end{tabular}\\[0.3cm]
\begin{tabular}{|l||c|c||c|c|}\hline
$h$& Godunov,$||c-c_h||_{L^1}$&$\alpha$&DFLU,$||c-c_h||_{L^1}$& $\alpha$\\
\hline
1/50&4.8407 $\times 10^{-2}$&& 4.8486 $\times 10^{-2}$&\\
\hline
1/100&3.0161 $\times 10^{-2}$&0.6825& 3.0201 $\times 10^{-2}$&0.6829\\
\hline
1/200& 1.9307 $\times 10^{-2}$&0.6435&1.9328$\times 10^{-2}$ &0.6439  \\
\hline
1/400&1.2618 $\times 10^{-2}$&0.6136&1.2628 $\times 10^{-2}$ &0.6140  \\
\hline
1/800&8.4125$\times 10^{-3}$&0.5848&8.4173 $\times 10^{-3}$&0.5851   \\
\hline
\end{tabular}
\end{center}
\caption{Riemann problem (\ref{flux1}), (\ref{IC2}): $L^1$-errors between exact and calculated solutions at $t=.5$.}
\label{table2}
\end{table}

\subsection{Comparison of the DFU, upstream mobility, FORCE and Lax-Friedrichs fluxes}
   In the previous section, we have seen that Godunov and DFLU fluxes give schemes with very close performances. In this section we compare the DFLU flux with the other fluxes that we mentioned in section \ref{finite difference} which are the
upstream mobility, FORCE and Lax-Friedrichs fluxes. We take now 
\be \begin{array}{l} 
f(s,c)=  \var_1 = \dfrac{\lambda_1(s,c)}{\lambda_1(s,c) + \lambda_2(s,c)} [ \var + (g_1-g_2)\lambda_2(s,c) ],\\
\lambda_1(s,c)=\dfrac{s^2}{.5+c}, \lambda_2(s,c)=(1-s)^2, \,\,g_1=2, g_1=1,  \varphi=0,\\
a(c)=.25c. 
\end{array} \label{f2} \ee
In all following experiments the discretization is such that $\Delta t=1/125$ and $h=1/100$.

\begin{figure}[H]
 \includegraphics[width=7.2cm]{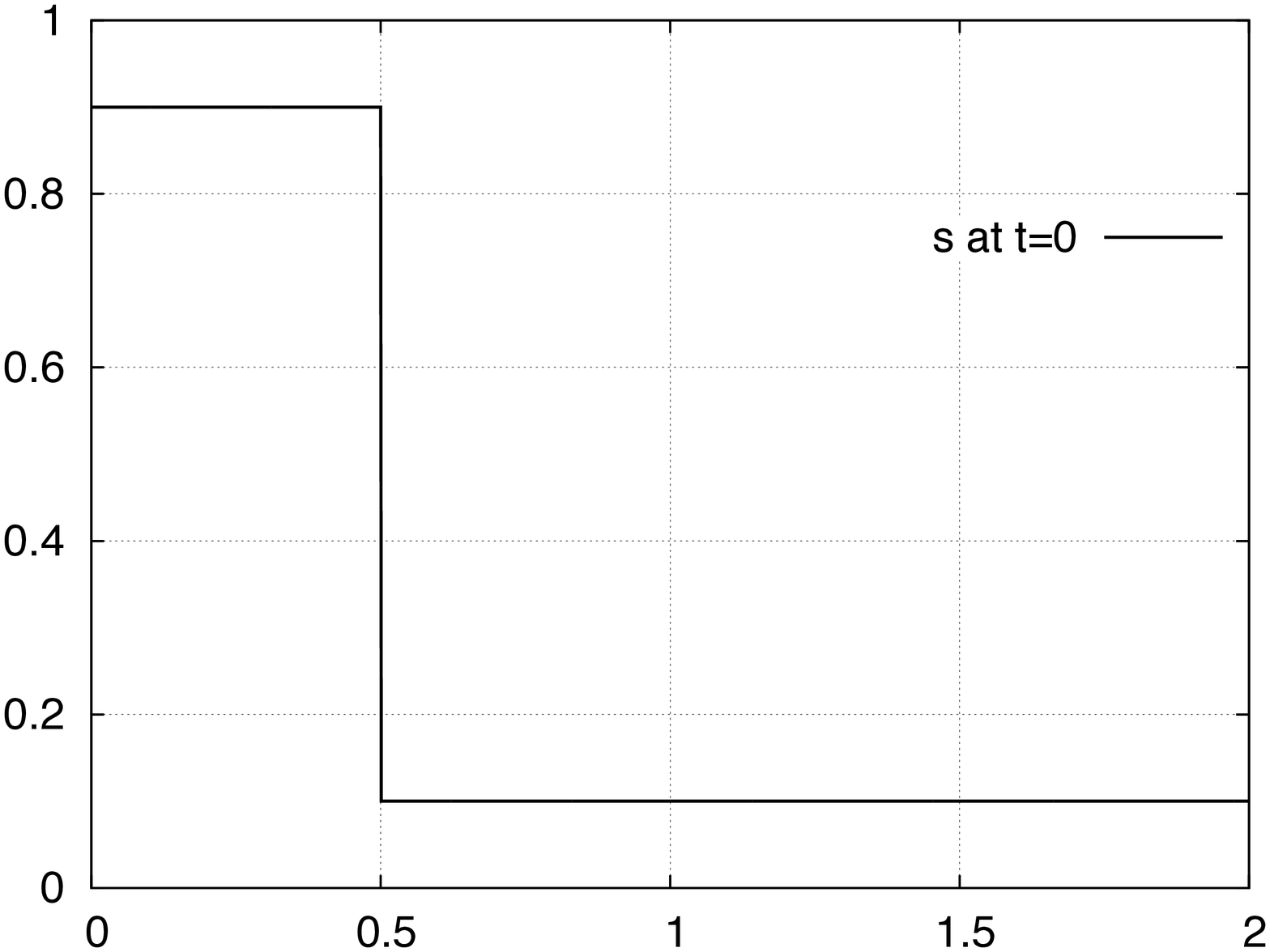}
 \includegraphics[width=7.2cm]{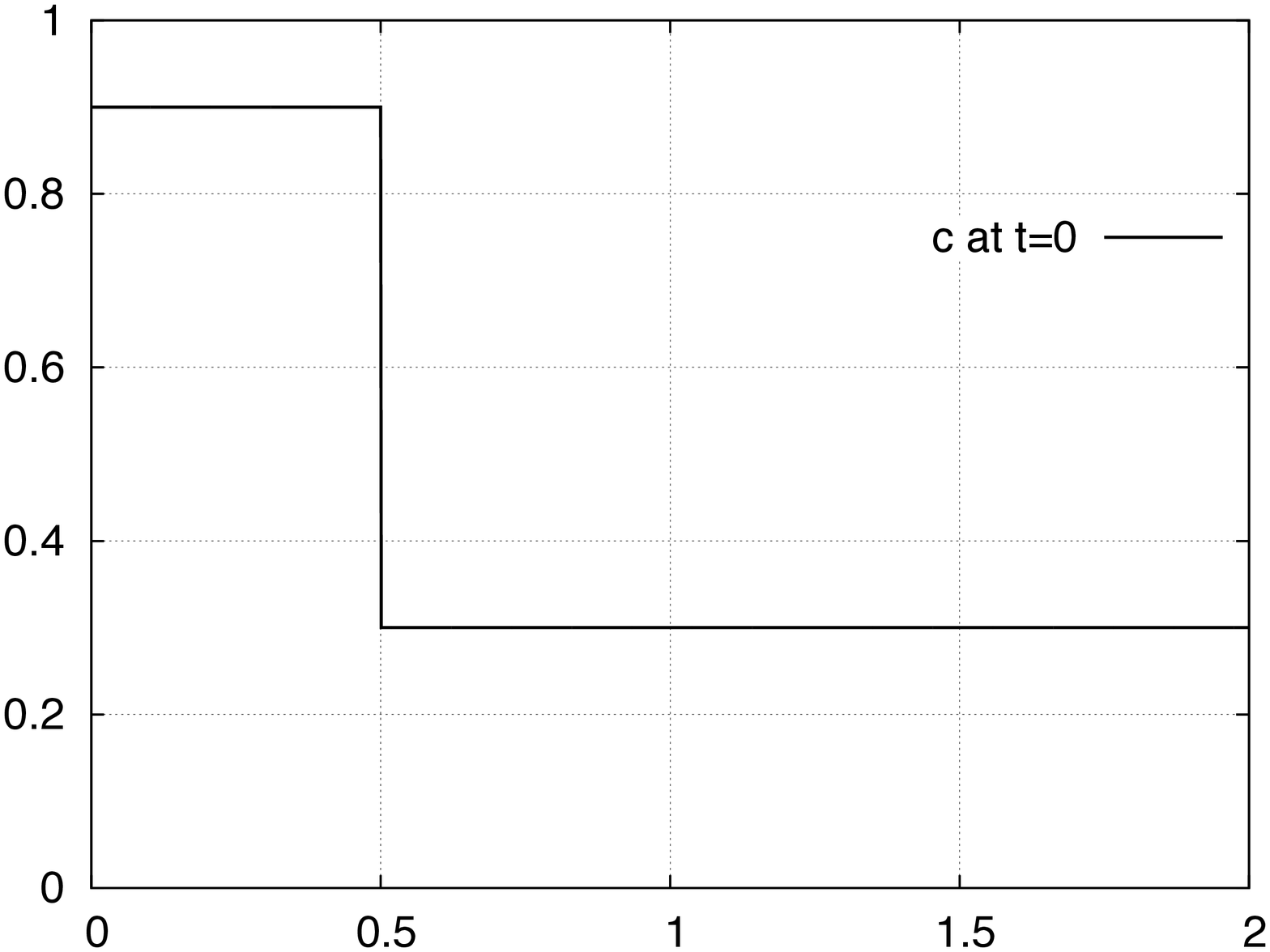}\\
 \includegraphics[width=7.2cm]{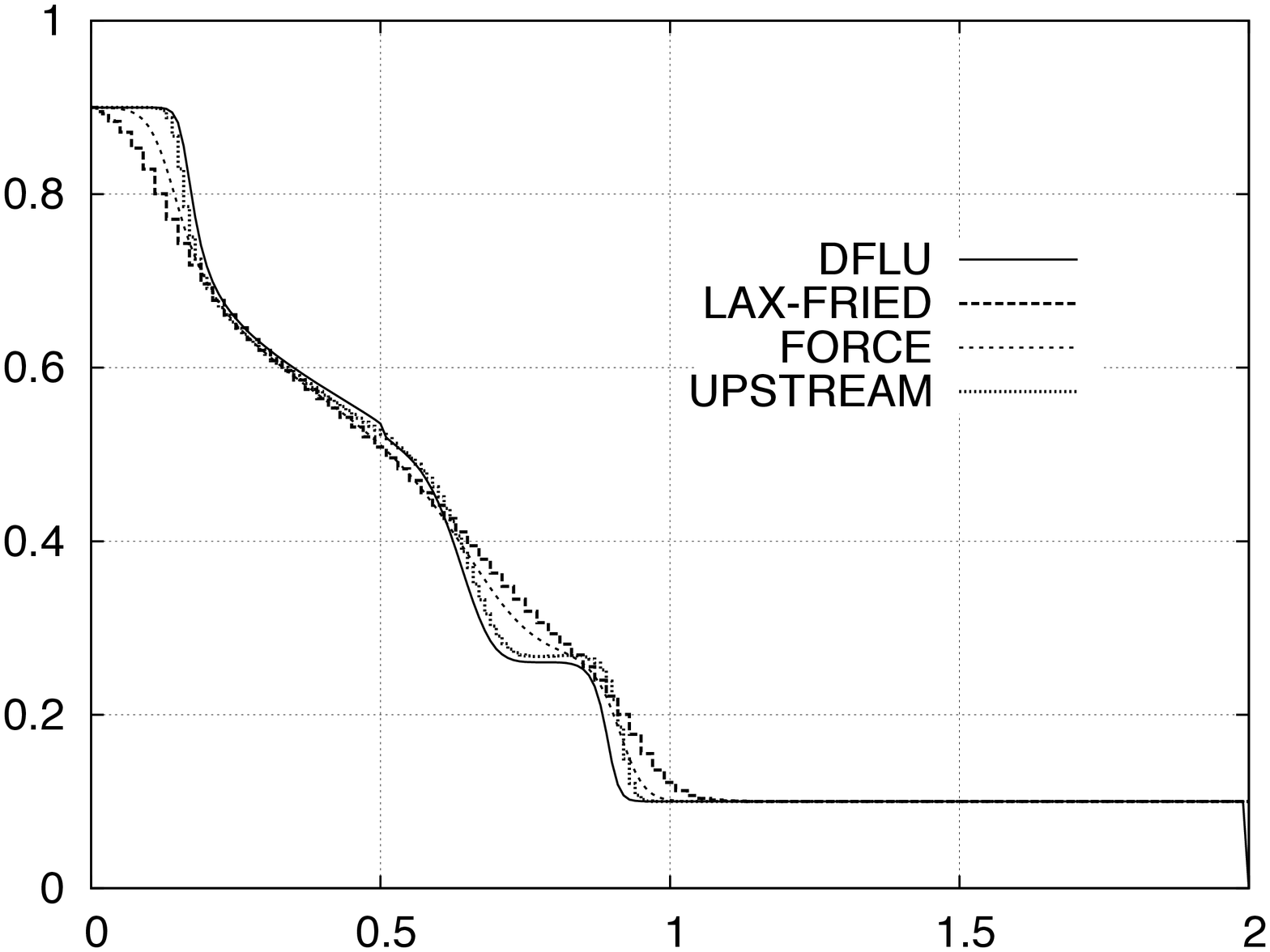}
 \includegraphics[width=7.2cm]{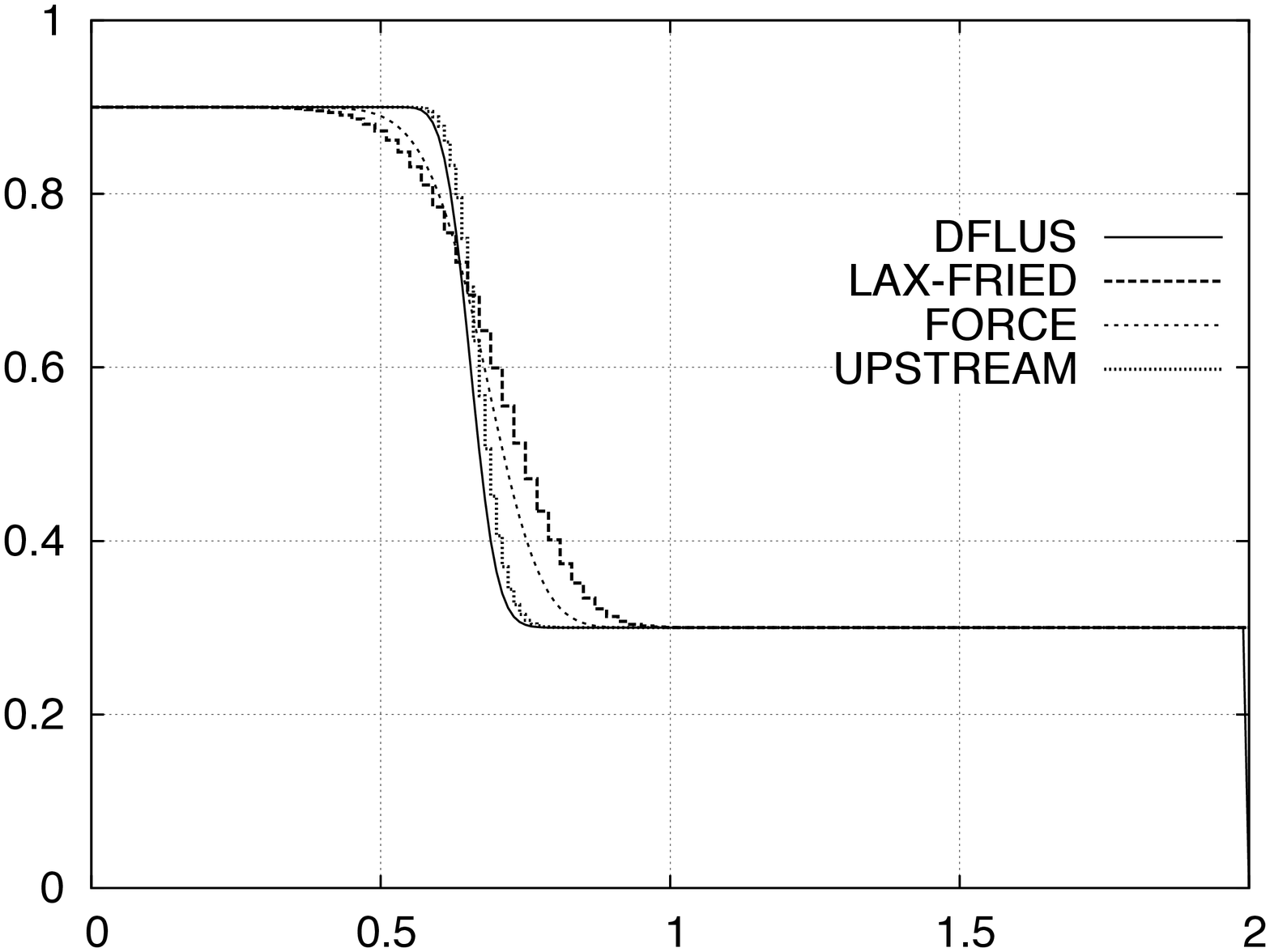}\\
 \includegraphics[width=7.2cm]{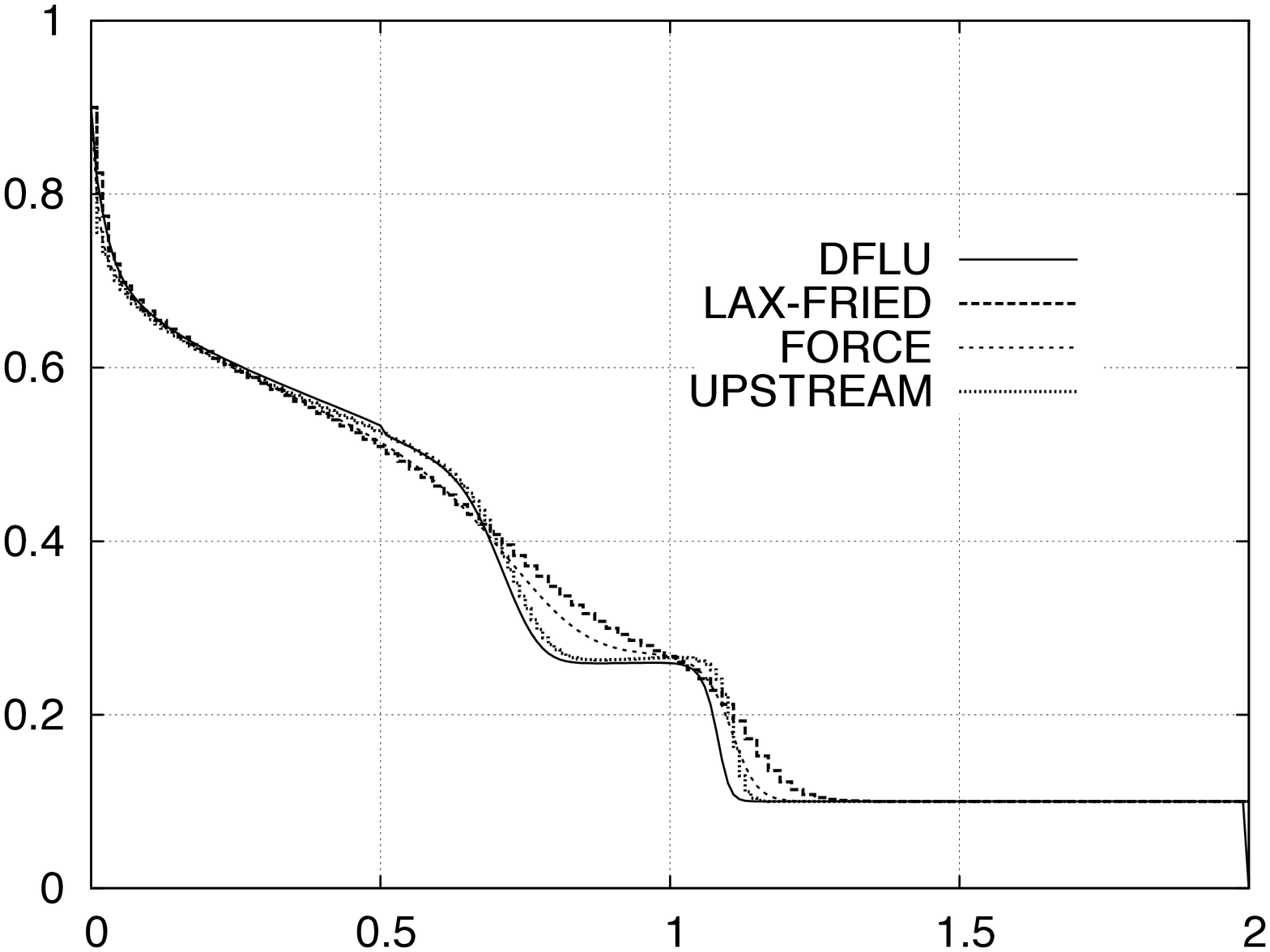}
 \includegraphics[width=7.2cm]{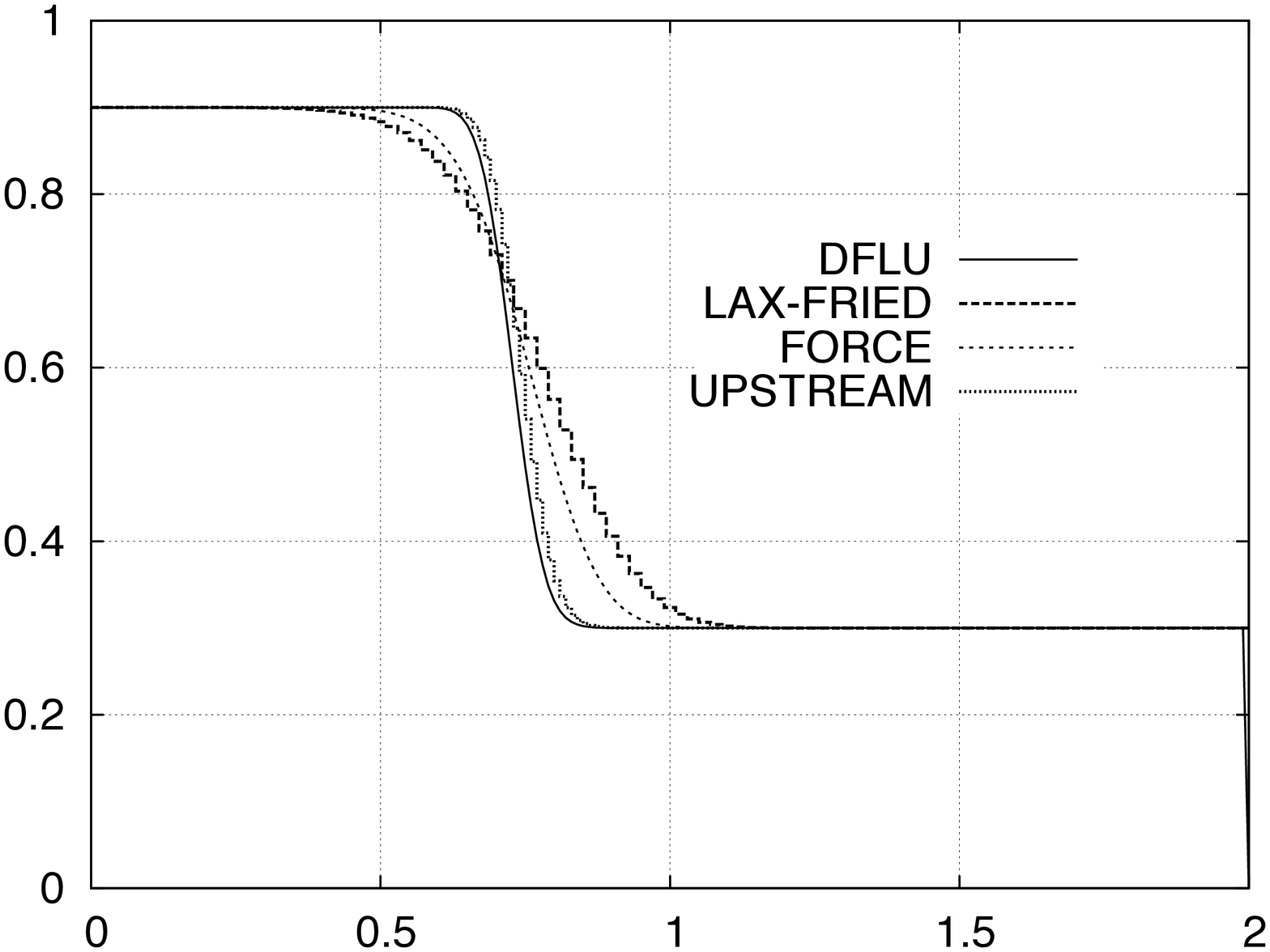}\\
\caption{$s$ (left) and $c$ (right) calculated at t=0., t=1. and t=1.5 for data (\ref{f2}), (\ref{ivpb}), (\ref{bvdiv}).}
\label{fig-154ivp}
\end{figure}

We first consider a pure initial value problem.
Initial condition (see top of Fig. \ref{fig-154ivp}) is given by   
\be s(x,0) = \left\{ \begin{array}{lll}
.9 &\mbox{if}& x<.5, \\ .1 &\mbox{if}& x>.5 \end{array} \right.  , \quad
c(x,0) = \left\{ \begin{array}{lll}
.9 &\mbox{if}& x<.5, \\ .3 &\mbox{if}& x>.5 \end{array} \right.  .
\label{ivpb}
\ee
With this initial condition we have $F^{DFLU}(s_L,s_R,c_L,c_R)=F^{G}(s_L,s_R,c_R,c_L)$ with
$ s_L=.9,s_R=.1,c_L=1.$ and $c_R=0$. Boundary data are such that
\be s(0,t)  = .9, \; s(2,t) = .1, \quad c(0,t)  = .9, \; c(2,t) = .3 \,\,\,\,\forall\,\,t \geq 0.
\label{bvdiv} \ee
Calculated solutions at time levels t=1 and 1.5 are shown in Fig.\ref{fig-154ivp}. They show that, as expected, the DFLU flux, which is the closest to a Godunov scheme, performs better than the other schemes. The upstream mobility flux, which is an upwind scheme, performs better than the  two central difference schemes, the FORCE and Lax-Friedrichs schemes. 

To confirm these first observations we consider now a boundary value problem. We just changed the boundary functions, so instead of boundary conditions (\ref{ivpb}) we consider now a problem with closed boundaries, that is fluxes are zero at the boundary:
\be f \equiv 0 \,\mbox{ at } \,\,x=0 \mbox{ and }  x=2 \,\,\mbox{ for all } \, t \geq 0  .
\label{bvpb} \ee
They show that, as expected, the DFLU scheme, which is the closest to a Godunov scheme, performs better than the  upstream mobility, the FORCE or the Lax-Friedrichs schemes. 

\begin{figure}[H]
 \includegraphics[width=7.2cm]{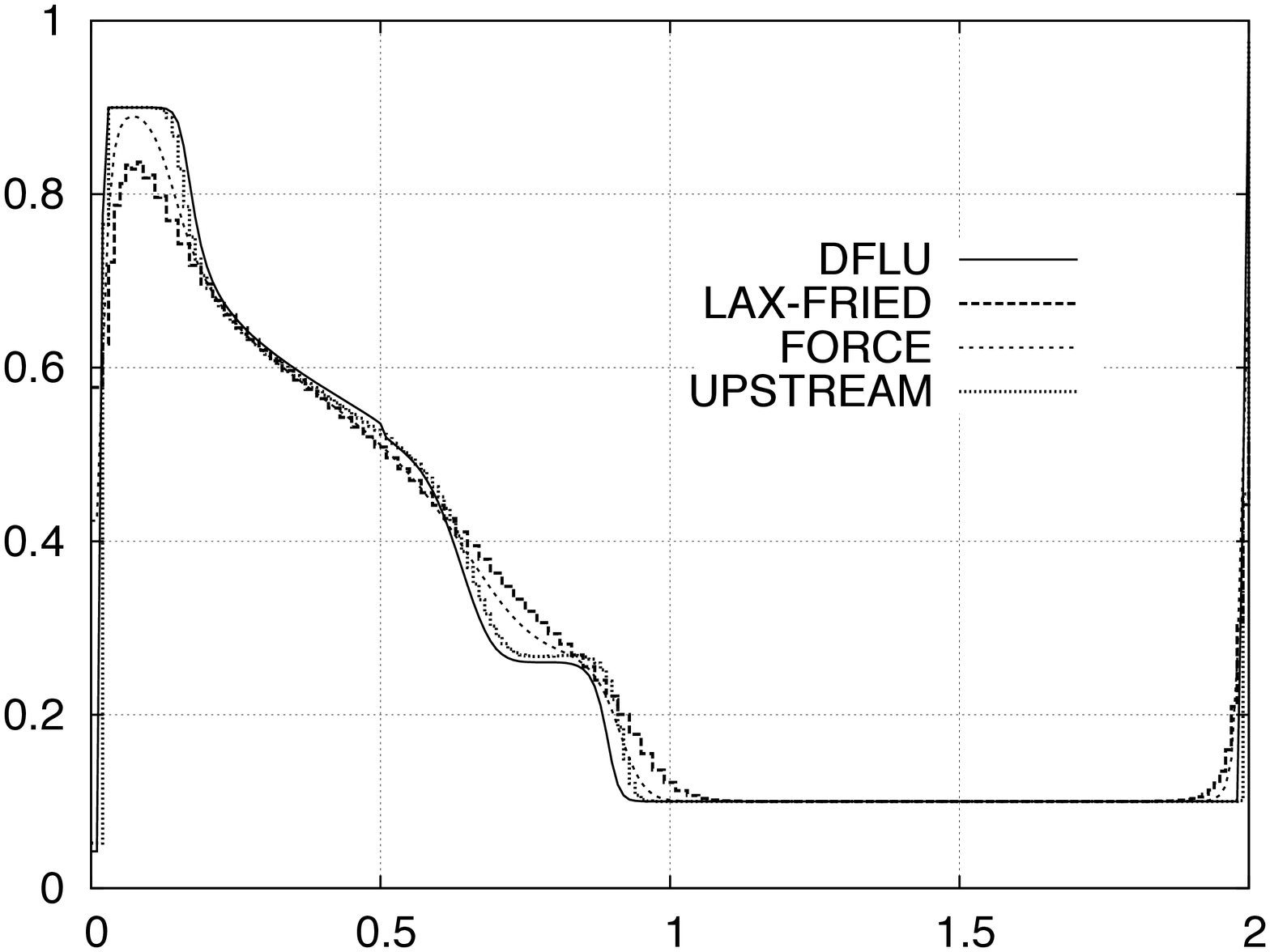}
 \includegraphics[width=7.2cm]{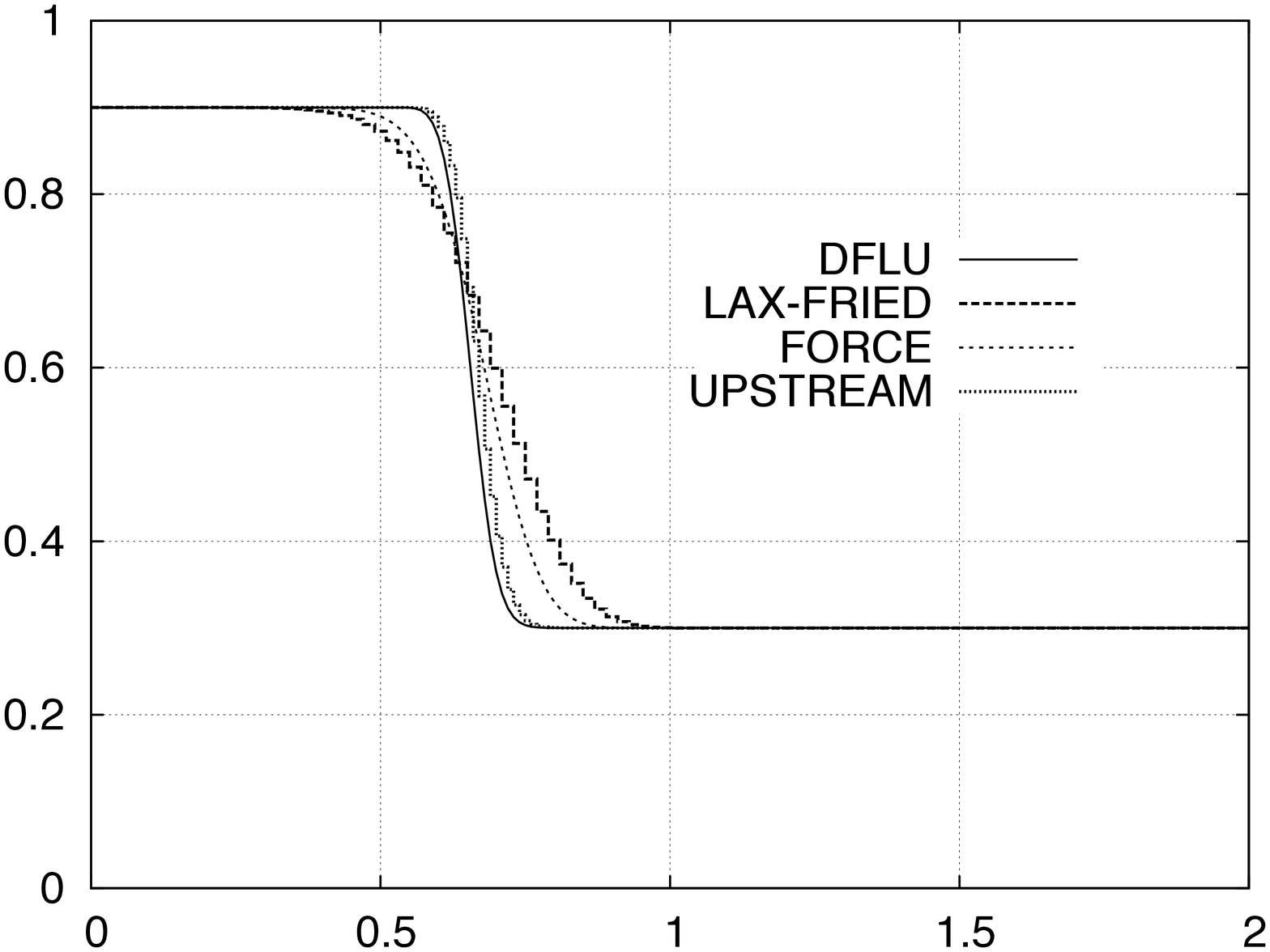}\\
 \includegraphics[width=7.2cm]{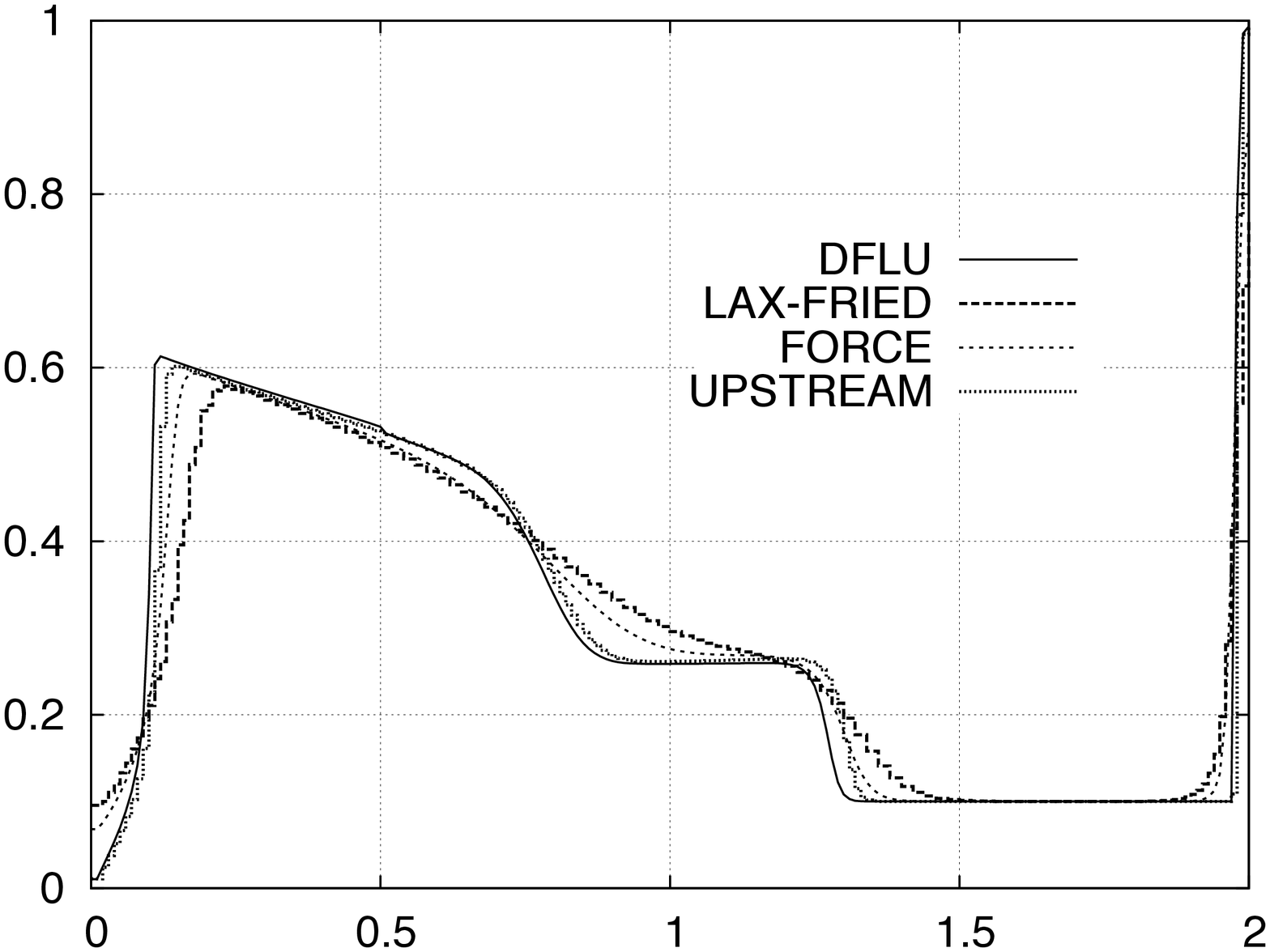}
 \includegraphics[width=7.2cm]{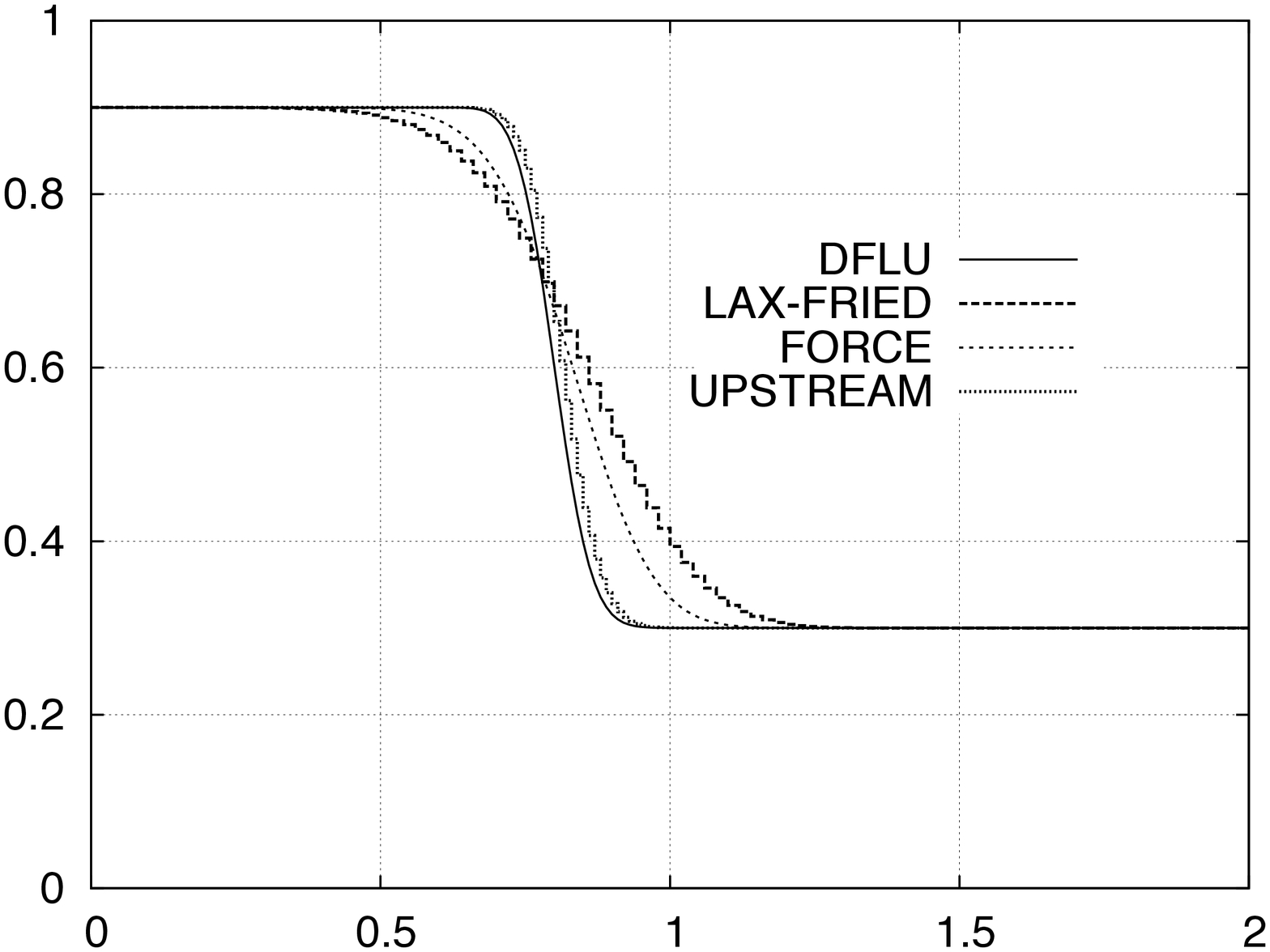}\\
 \includegraphics[width=7.2cm]{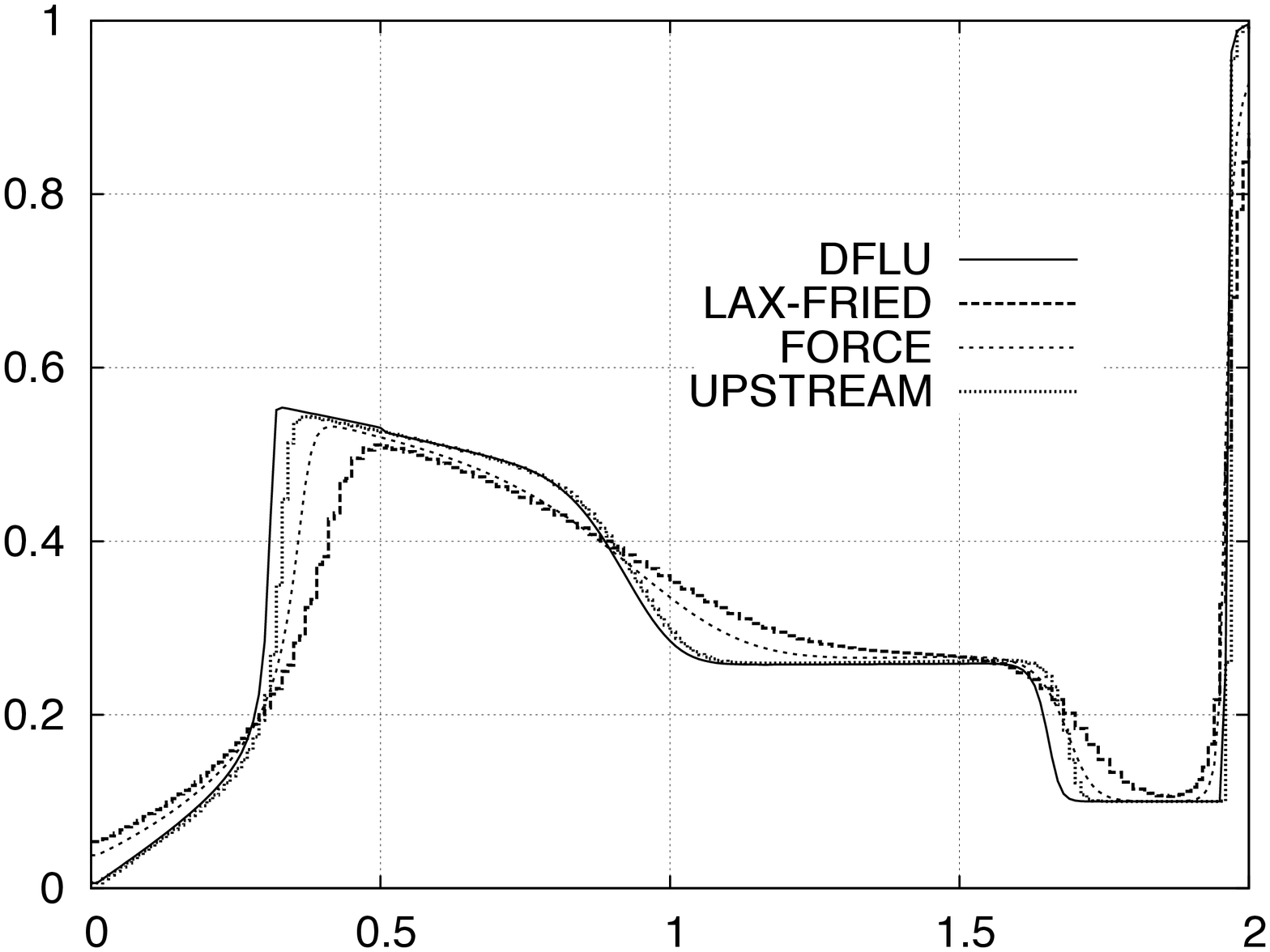}
 \includegraphics[width=7.2cm]{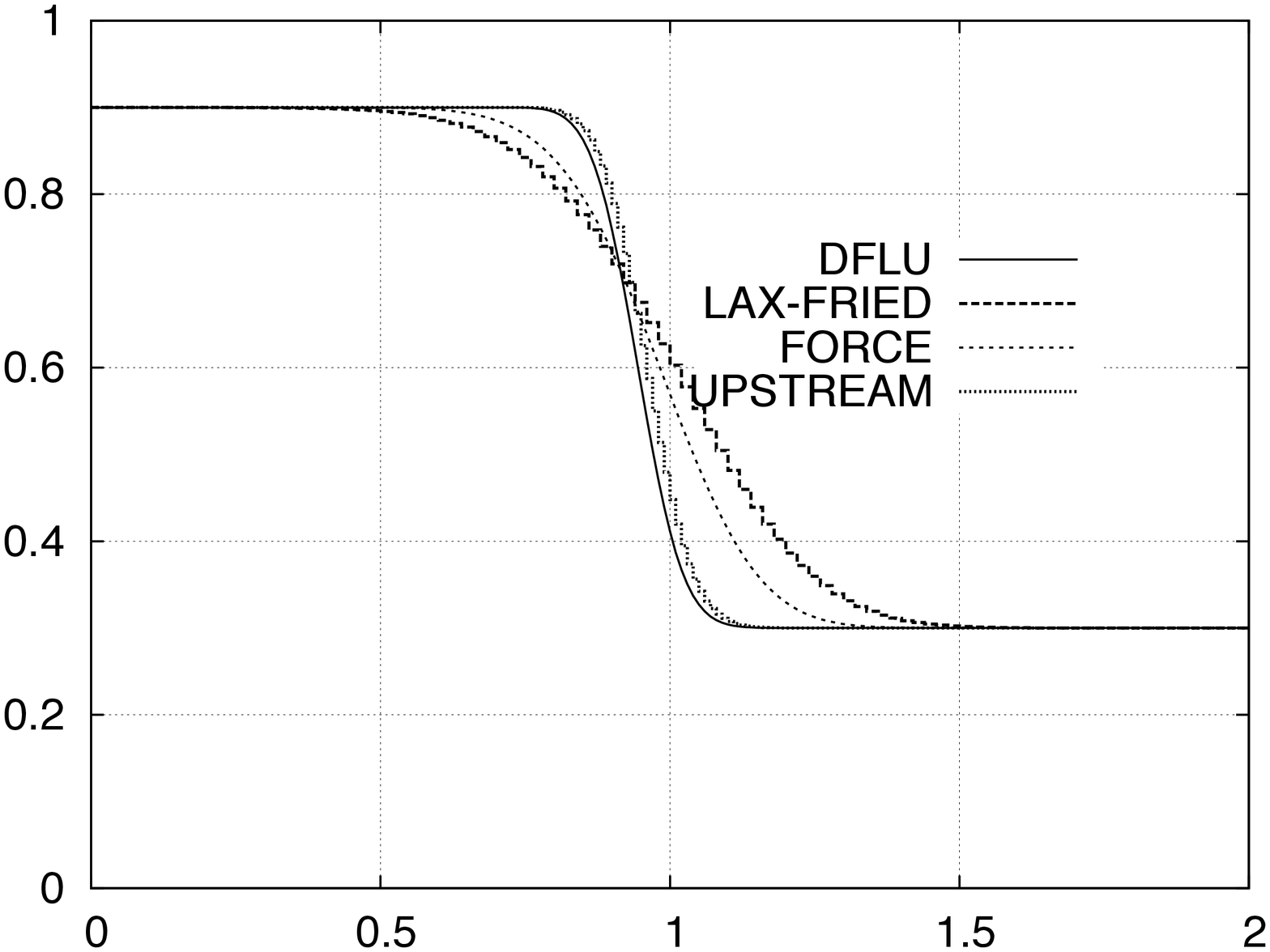}
 \caption{$s$ (left) and $c$ (right) calculated at t=1., t=2. and t=3. for data (\ref{f2}), (\ref{bvpb}), (\ref{bvdiv}).}
\label{fig-bvp}
\end{figure}

\begin{figure}[H]
 \includegraphics[width=7.2cm]{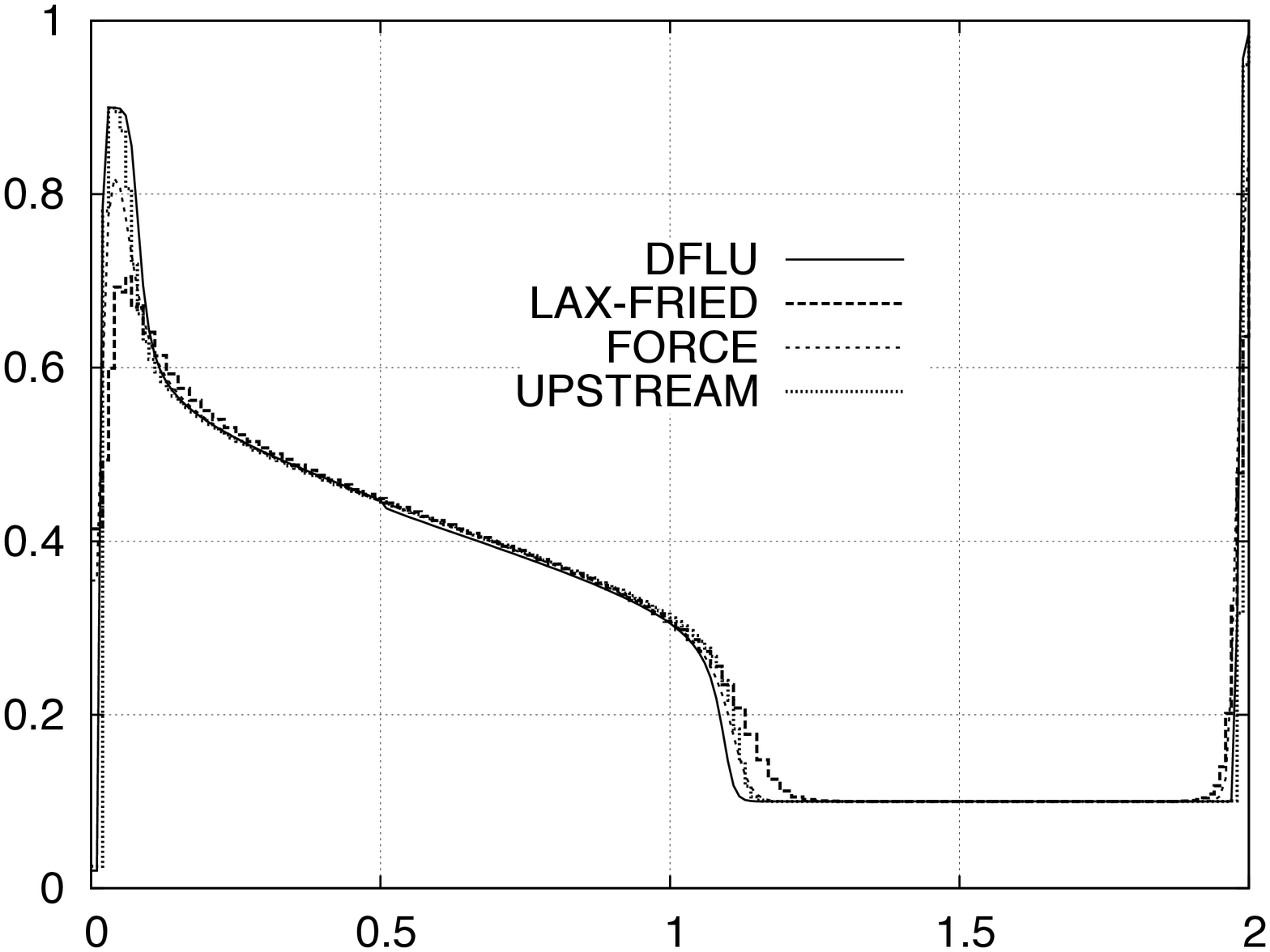}
 \includegraphics[width=7.2cm]{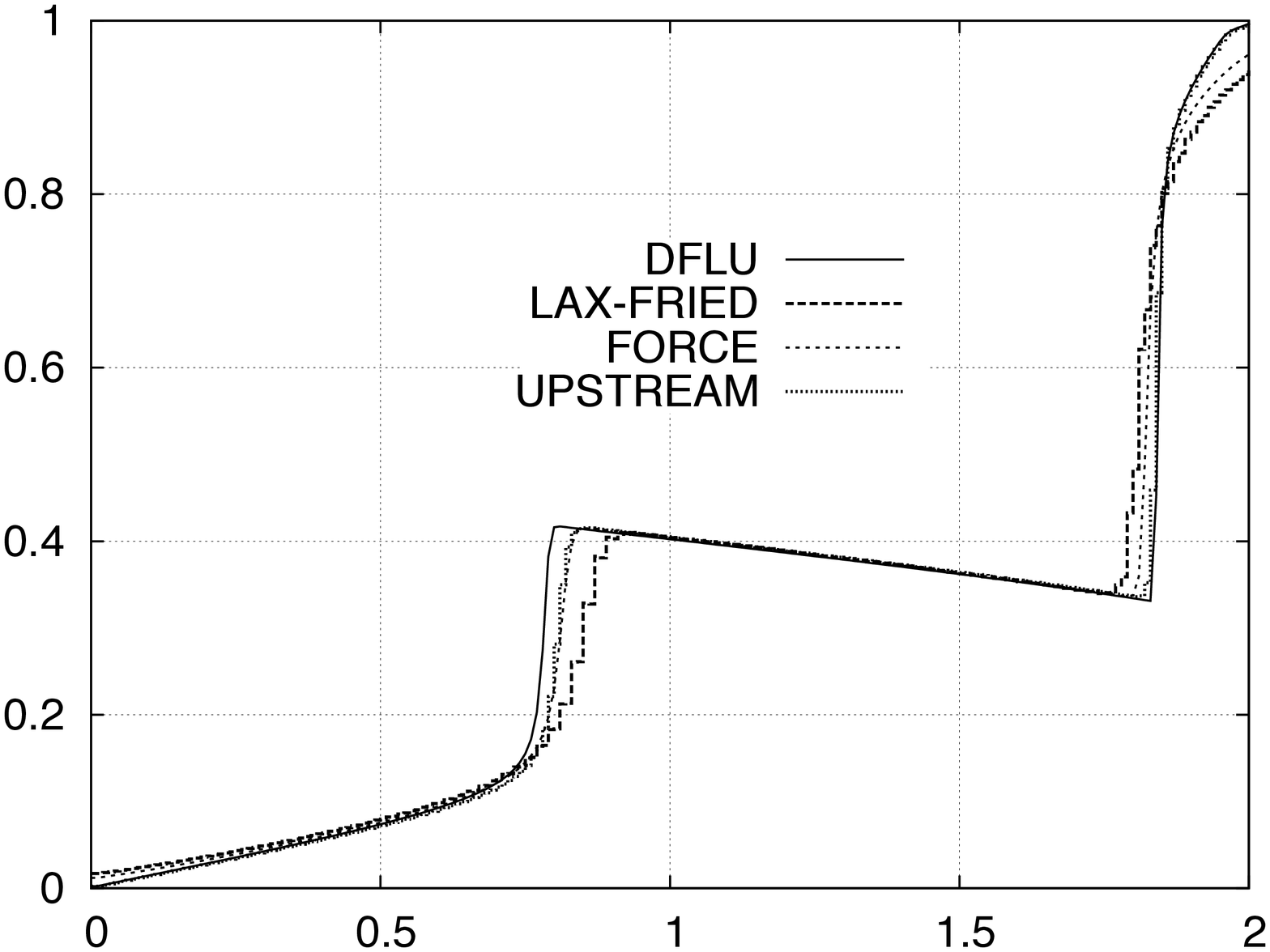}
\caption{$s$ (left) calculated at t=1. and t=3. for same data as in Fig. \ref{fig-bvp} but without polymer injection.}
\label{figbvpc0u}
\end{figure}

The purpose of the last experiment whose results are shown in Fig. \ref{figbvpc0u} is to show the effect of polymer flooding. In this experiment we remove polymer flooding  and take  $c\equiv 0$ at all time. By comparing with the solution shown in Fig. \ref{fig-bvp} bottom left we observe that as expected the saturation front is moving faster since there is no retardation due to the increase of viscosity of the wetting fluid caused by the polymer injection. We also observe that the structure of the solution is less complex.

\section{Conclusion}
The DFLU flux defined in \cite{AdiJafGow04} for scalar conservation laws was used to construct a new scheme for a class of system of conservation laws. It was applied to a system for polymer flooding. It is very close to the flux given by an exact Riemann solver and the corresponding finite volume scheme compares favorably to other schemes using the uptream mobility, the Lax-Friedrichs and the FORCE fluxes. The DFLU is also very easy to implement. The extension to the case with a change of rock type is straightforward since the DFLU flux was built to solve this case. It will work even in cases where the upstream mobility fails \cite{SidJaf09}.
In a separate paper \cite{AdiJafGow09c} we show how to use the DFLU flux to solve Hamilton-Jacobi equations with a discontinuous Hamiltonian.
\bibliographystyle{siam}
  
\bibliography{hyper}

\end{document}